\documentclass[aap,seceqn,MSNbibl,citesort,dvips]{arximspdf}

\doi{10.1214/10-AAP716}
\volume{21}
\issue{3}
\pubyear{2011}
\firstpage{1017}
\lastpage{1052}

\makeatletter

\newtheorem{theo}{Theorem}[section]
\newtheorem{lem}[theo]{Lemma}
\newtheorem{cor}[theo]{Corollary}

\newproclaim{rem}[theo]{Remark}
\newproclaim{defi}[theo]{Definition}

\newcommand{\EE}{\mathbb{E}}
\newcommand{\LL}{\mathbb{L}}
\newcommand{\RR}{\mathbb{R}}
\newcommand{\Ba}{\mathcal{B}}
\newcommand{\Da}{\mathcal{D}}
\newcommand{\Ea}{\mathcal{E}}
\newcommand{\Sa}{\mathcal{S}}
\newcommand{\Ra}{\mathcal{R}}
\newcommand{\Va}{\mathcal{V}}
\newcommand{\Fa}{\mathcal{F}}
\newcommand{\Ga}{\mathcal{G}}
\newcommand{\Xa}{\mathcal{X}}
\newcommand{\Ma}{\mathcal{M}}
\newcommand{\Pa}{\mathcal{P}}
\newcommand{\point}{\cdot}
\newcommand{\et}{\eta}

\makeatother

\begin{document}
\begin{frontmatter}

\title{Concentration inequalities for mean field particle~models}
\runtitle{$\!\!\!$Concentration inequalities for mean field particle models}

\begin{aug}
\author[A]{\fnms{Pierre} \snm{Del Moral}\corref{}\ead[label=e1]{Pierre.Del-Moral@inria.fr}} and
\author[B]{\fnms{Emmanuel} \snm{Rio}\ead[label=e2]{rio@math.uvsq.fr}}
\runauthor{P. Del Moral and E. Rio}
\affiliation{INRIA Bordeaux-Sud-Ouest, and INRIA Bordeaux-Sud-Ouest and~Universit\'{e}~de~Versailles}
\address[A]{Centre INRIA Bordeaux-Sud-Ouest\\
Institut de Math\'{e}matiques de Bordeaux\\
Universit\'{e} Bordeaux 1\\
351 Cours de la Lib\'{e}ration\\
33405 Talence Cedex\\
France\\
\printead{e1}} 
\address[B]{Centre INRIA Bordeaux-Sud-Ouest\\
and\\
Laboratoire de Math\'{e}matiques de Versailles\\
Universit\'{e} de Versailles, B\^{a}timent Fermat\\
45 Av. des Etats-Unis\\
78035 Versailles Cedex\\
France\\
\printead{e2}}
\end{aug}

\received{\smonth{4} \syear{2010}}

%
\begin{abstract}
This article is concerned with the fluctuations and the concentration
properties of a general class of discrete generation and mean field
particle interpretations of nonlinear measure valued processes. We
combine an original stochastic perturbation analysis with a
concentration analysis for triangular arrays of conditionally
independent random sequences, which may be of independent interest.
Under some additional stability properties of the limiting measure
valued processes, uniform concentration properties, with respect to the
time parameter, are also derived. The concentration inequalities
presented here generalize the classical Hoeffding, Bernstein and
Bennett inequalities for independent random sequences to interacting
particle systems, yielding very new results for this class of models.

We illustrate these results in the context of McKean--Vlasov-type
diffusion models, McKean collision-type models of gases and of a class
of Feynman--Kac distribution flows arising in stochastic engineering
sciences and in molecular chemistry.
\end{abstract}

%

\end{frontmatter}

\section{Introduction}\label{intro}

\subsection{Mean field particle models}

Let $(E_n)_{n\geq0}$ be a sequence of measurable spaces equipped with
some $\sigma$-fields $(\Ea_n)_{n\geq0}$, and we let $\Pa(E_n)$ be the
set of all probability measures over the set $E_n$,
with $n\geq0$. We consider a collection
of transformations $\Phi_{n+1}\dvtx\Pa(E_{n})\rightarrow\Pa(E_{n+1})$,
$n\geq0$, and we denote by $(\eta_n)_{n\geq0}$
a sequence of probability measures on $E_n$ satisfying a nonlinear
equation of the following form:
%
%
\begin{equation}\label{phi}
\eta_{n+1}=\Phi_{n+1}(\eta_{n}) .
\end{equation}
These discrete time versions of conservative and nonlinear
integro-differential type equations
in distribution spaces
arise in a variety of scientific disciplines including in physics,
biology, information theory and engineering sciences.
To motivate the article, before describing their mean field particle
interpretations,
we illustrate these rather abstract evolution models
working out explicitly some of these equations in
a series of concrete examples. The first one is related to nonlinear filtering
problems arising in signal processing. Suppose we are given
a pair signal-observation Markov chain $(X_n,Y_n)_{n\geq0}$ on some
product space $(\RR^{d_1}\times\RR^{d_2})$,
with some initial distribution and
Markov transition of the following form:
\begin{eqnarray*}
\mathbb{P}\bigl((X_0,Y_0)\in d(x,y)\bigr)&=&\eta_0(dx) g_0(x,y) \lambda
_0(dy),\\
\mathbb{P}\bigl((X_{n+1},Y_{n+1})\in d(x,y)|(X_n,Y_n)
\bigr)&=&M_{n+1}(X_n,dx) g_{n+1}(x,y) \lambda_{n+1}(dy).
\end{eqnarray*}
In the above display,
$\lambda_{n}$ stands for
some reference probability measures on~$\RR^{d_2}$, $g_n$ is a
sequence of positive functions, $M_{n+1}$ are
Markov transitions from $\RR^{d_1}$ into itself and finally $\eta_0$
stands for
some initial probability measure on $\RR^{d_1}$. For a given sequence of
observations $Y=y$ delivered by some sensor, the filtering problem
consists of computing sequentially
the flow of conditional distributions defined by
\[
\widehat{\eta}_{n}=\operatorname{Law}(X_{n}|Y_0=y_0,\ldots
,Y_{n}=y_{n})
\]
and
\[
\eta_{n+1}=\operatorname{Law}(X_{n+1}|Y_0=y_0,\ldots,Y_{n}=y_{n}).
\]
These distributions satisfy
a nonlinear evolution equation of the form~(\ref{phi})
with the transformations
%
%
\begin{eqnarray}\label{exref}
\Phi_{n+1}(\eta_{n})(dx')&=&\int\widehat{\eta
}_{n}(dx) M_{n+1}(x,dx') \quad\mbox{and}\nonumber\\[-8pt]\\[-8pt]
\widehat{\eta}_{n}(dx)&=&\frac{G_n(x)}{\int\eta_{n}(dx') G_n(x')}
\eta_{n}(dx)\nonumber
\end{eqnarray}
for some collection of likelihood functions $G_n=g_n(\cdot,y_n)$.
Replacing these functions by some $]0,1]$-valued
potential function $G_n$ on $\RR^{d_1}$, we obtain the conditional
distributions of a
particle absorption model $X'_n$ with free evolution transitions $M_n$
and killing rate $(1-G_n)$. More precisely, if $T$ stands for the
killing time of the process,
we have that
%
%
\begin{equation}\label{exrefab}
\widehat{\eta}_{n}=\operatorname{Law}(X_{n}' | T> n) \quad
\mbox{and}\quad
\eta_{n+1}=\operatorname{Law}(X_{n+1}' | T>n).
\end{equation}
These nonabsorption conditional distributions arise in the analysis of
confinement processes, as well as in
computational physics with the numerical solving of Schr\"{o}dinger
ground state energies (see, e.g., \cite{Miclo3,DDS,Rousset2}).

Another important class of measures arising in particle physics and
stochastic optimization problems
is the class of Boltzmann distributions,\vadjust{\goodbreak} also known as the Gibbs
measure, defined by
%
%
\begin{equation}\label{exrefabg}
\eta_n(dx)=\frac{1}{Z_n} e^{-\beta_n V(x)} \lambda(dx) ,
\end{equation}
with some reference probability measure $\lambda$, some inverse
temperature parameter
and some nonnegative potential energy function $V$ on some state
space $E$. The normalizing constant
$Z_n$ is sometimes called ``partition function'' or the ``free
energy.''
In sequential Monte Carlo methodology,
as well as in operation research literature, these
multiplicative formulae are often used to compute rare events
probabilities, as well as the cardinality
or the volume of some complex state spaces. Further details on these stochastic
techniques can be found in the series of articles \mbox{\cite{cerou,Rubin1,Rubin2}}.
To fix the ideas we can consider the uniform
measure on some finite set $E$. Surprisingly, this flow of measures
also satisfies the above nonlinear
evolution equation, as soon as $\eta_n$ is an invariant measure of
$M_n$, for each $n\geq1$, and the potential functions
$G_n$ are chosen of the following form:
$G_{n}=\exp{((\beta_{n+1}-\beta_n)V)}$. The partition
functions can also be computed in terms
of the flow of measures $(\eta_p)_{0\leq p<n}$ using the easy-to-check
multiplicative formula,
\[
Z_n=\prod_{0\leq p<n}\int\eta_p(dx) G_p(x),
\]
as soon as $\beta_0=0$. In statistical mechanics literature, the above formula
is sometimes called the Jarzynski or the Crooks equality \cite
{Cro,Jar,Jar2}. Notice that the stochastic models
discussed above remain valid if we replace $ e^{-\beta_{n+1} V}$ by any
collection of functions $g_{n+1}$ s.t.
$
g_{n+1}=g_{n} \times G_{n}=\prod_{0\leq p\leq n}G_p
$,
for some potential functions $G_n$, with $n\geq0$. Further details on this
model, with several worked-out applications on concrete hidden Markov
chain problems and Bayesian inference
can be found in the article \cite{ddj} dedicated to sequential Monte
Carlo technology. All of the models discussed above can be abstracted
in a single probabilistic
model. The latter is often called a Feynman--Kac model, and it will be
presented in some details in Section~\ref{secFKmod}.

The mean field-type interacting particle system associated with
equation~(\ref{phi})
relies on the fact that the one-step mappings can be rewritten in the
following form:
%
%
\begin{equation}\label{phidef}
\Phi_{n+1}(\eta_{n})=\eta_{n}K_{n+1,\eta_{n}}
\end{equation}
for some collection of Markov kernels $K_{n+1,\mu}$ indexed by the time
parameter $n\geq0$, and the set of measures $\mu$ on the space
$E_{n}$. We already mention that the choice of the Markov transitions
$K_{n,\eta}$ is not unique.
Several examples are presented in Section~\ref{secillus} in the context
of Feynman--Kac semigroups
or McKean--Vlasov diffusion-type models. To fix the ideas, we can choose
elementary transitions $K_{n,\eta}(x,dx')=\Phi_{n}(\eta
)(dx')$ that do not depend on the state variable $x$.

In the literature on mean field particle models, the transitions
$K_{n,\eta}$ are called a choice of McKean transitions.
These models provide a natural interpretation
of the flow of measures $\eta_n$ as the laws of the time inhomogeneous
Markov chain
$\overline{X}_n$ with elementary transitions
\[
\mathbb{P}(\overline{X}_n\in dx|\overline{X}_{n-1})=K_{n,\eta
_{n-1}}(\overline{X}_{n-1},dx) \qquad\mbox{with }
\eta_{n-1}=\operatorname{Law}(\overline{X}_{n-1})
\]
and starting with some initial random variable with
distribution
$\eta_0=\operatorname{Law}(\overline{X}_0)$. The Markov chain $\overline
{X}_n$ can be thought of as a perfect
sampling algorithm. For a thorough description of these discrete
generation and nonlinear McKean-type models, we refer the reader
to \cite{fk}.
In the further development of the article, we always assume that
the mappings
\[
(x_n^i)_{1\leq i\leq N}\in E^N_n\mapsto K_{n+1,
{1/N}\sum_{j=1}^N\delta_{x^j_n}}(x^i_n,A_{n+1})
\]
are $\Ea^{\otimes N}_n$-measurable, for any $n\geq0$, $N\geq1$, and
$1\leq i\leq N$, and any measurable subset $A_{n+1}\subset E_{n+1}$. In
this situation, the mean field particle interpretation of this
nonlinear measure valued model is an $E^N_n$-valued Markov chain
$\xi^{(N)}_n=(\xi^{(N,i)}_n)_{1\leq i\leq N}$,
with elementary transitions
defined as
%
%
\begin{eqnarray}\label{meanfield}
\mathbb{P}\bigl(\xi^{(N)}_{n+1}\in dx
| \Fa^{(N)}_n\bigr) =\prod_{i=1}^N K_{n+1,\eta^N_n}\bigl(
\xi^{(N,i)}_n,dx^i\bigr) \nonumber\\[-8pt]\\[-8pt]
&&\eqntext{\mbox{with } \eta^N_n:=\displaystyle \frac{1}{N}
\sum_{j=1}^N\delta_{\xi_n^{(N,j)}}.}
\end{eqnarray}
In the above displayed formula, $\Fa^{N}_n$ stands for the $\sigma$-field
generated by the random sequence $(\xi^{(N)}_p)_{0\leq p\leq n}$, and
$dx=dx^1\times\cdots\times dx^N$ stands
for an infinitesimal neighborhood of a point $x=(x^1,\ldots,x^N)\in E_n^N$.
The initial system $\xi^{(N)}_0$ consists of $N$ independent and
identically
distributed random variables with common law~$\eta_0$.
As usual, to simplify the presentation, when there is no possible
confusion we suppress the
parameter $N$, so that we write $\xi_n$ and $\xi^i_n$ instead of $\xi
^{(N)}_n$ and $\xi^{(N,i)}_n$.
The state components of this
Markov chain are called particles or sometimes walkers in physics to
distinguish the
stochastic sampling model with the physical particle in molecular models.

The rationale behind this is that $\eta^N_{n+1}$ is the empirical
measure associated
with $N$ independent variables with distributions $K_{n+1,\eta
^N_n}
(\xi^{i}_n,dx)$, so
as soon as $\eta^N_{n}$ is a good approximation of $\eta_{n}$ then, in
view of
(\ref{meanfield}), $\eta^N_{n+1}$ should be a good approximation of
$\eta_{n+1}$. Roughly speaking, this induction
argument shows that $\eta^N_n$ tends to $\eta_n$, as the population
size $N$ tends to infinity.

These stochastic particle algorithms can be thought of
in various ways:
from the physical view point, they can be seen
as microscopic particle
interpretations of physical nonlinear measure valued equations.
From the pure
mathematical point of view, they can also be interpreted as natural stochastic
linearizations of nonlinear evolution semigroups. From the
probabilistic point of view, they can be interpreted as
interacting recycling acceptance-rejection sampling techniques.
In this case, they can be seen as a sequential and interacting importance
sampling technique.

For instance, in the context of the nonlinear filtering equation~(\ref{exref}),
the mean field particle model associated with the flow of optimal
one-step predictors $\eta_n$,
with the McKean transitions $K_{n,\eta}(x,dx')=\Phi_{n}(\eta
)(dx')$,
is the $(\RR^{d_1})^N$-valued Markov chain
defined by sampling $N$ conditionally independent random variables
$
\xi_{n+1}=(\xi^i_{n+1})_{1\leq i\leq N}
$, with common distribution given by
%
%
\begin{equation}\label{meanfield0}
\Phi_{n+1}\Biggl(\frac{1}{N}
\sum_{j=1}^N\delta_{\xi_n^{j}}\Biggr)(dx)=\sum_{i=1}^N\frac
{G_n(\xi
^i_n)}{\sum_{j=1}^NG_n(\xi^j_n)} M_{n+1}(\xi^i_n,dx).
\end{equation}
By construction, the resulting particle model is a simple genetic-type
stochastic algorithm:
the mutation and the selection transitions are dictated by the
prediction and the updating transitions defined in~(\ref{exref}).
During the
selection transition, one updates the positions of the particles in
accordance with the fitness likelihood functions
$G_n$.
This mechanism is called the selection-updating transition as the more
likely particles with high $G_n$-potential value
are selected
for reproduction. In other
words, this transition allows particles to give birth to some particles
at the expense
of light particles which die.
The second mechanism is called the mutation-prediction transition since
at this step each particle
evolves randomly according to the transition kernels $M_n$.
Another important feature of genetic-type particle models is that their
ancestral or their complete genealogical tree structure
can be used
to approximate the smoothing problem, including the computation of the
distribution of the signal trajectories given the observations. Further
details on this subject can be found in \cite{fk,DDS}.

The same genetic-type
particle algorithm applies for the particle absorption model (\ref
{exrefab}) and the Boltzmann--Gibbs
model~(\ref{exrefabg}), by replacing, respectively, the likelihood
functions by the nonabsorption rates
$G_n$ and the fitness functions $G_{n}=\exp{((\beta_{n+1}-\beta
_n)V)}$.

In the reverse angle, the occupation
measures of a given genetic-type particle mean field model converge, as
the size of the population tends to infinity, to the solution of an
evolution equation
of the form~(\ref{phi}), with the one-step transformations (\ref
{exref}). These limiting models are often called the infinite
population models. For a recent treatment on these genetic models, we
refer the reader to \cite{steins}.

The origins of genetic-type particle methods can be traced back in
physics and molecular chemistry
in the 1950s with the pioneering works of Harris and Kahn \cite
{harris} and Rosenbluth and Rosenbluth \cite{rosen}. During the last two
decades, the mean field particle interpretations of these discrete generation
measure valued equations were increasingly identified as a powerful
stochastic simulation algorithm
with emerging subjects in physics, biology and engineering sciences.
They have led to spectacular
results in signal processing processing with the corresponding particle
filter technology,
in stochastic engineering with interacting-type Metropolis and Gibbs
sampler methods,
as well as in quantum chemistry with quantum and diffusion Monte Carlo
algorithms
leading to precise estimates of the top eigenvalues
and the ground states of Schroedinger operators. For a thorough
discussion on these application areas,
we refer the reader to \cite{fk,ddj,arnaud}, and the references
therein. To motivate the article, we illustrate
the fluctuation and the concentration results presented
in this work with three additional illustrative examples, including
Feynman--Kac models,
McKean--Vlasov diffusion-type models, as well as interacting jump type
McKean model of gases.

We end this Introduction with some more or less traditional notation
used in the present article. We denote, respectively, by
$\mathcal{M}(E)$, $\mathcal{M}_{0}(E)$ and $\mathcal{B}(E)$, the set of
all finite signed measures on some measurable space $(E,\mathcal{E})$,
the convex subset of measures with null mass and the Banach space of
all bounded and measurable functions $f$ equipped with the uniform norm
$\Vert f\Vert$. We also denote by $\operatorname{Osc}_{1}(E)$, the convex set
of $\mathcal{E}$-measurable functions $f$ with oscillations
$\mbox{osc}(f)\leq1$. We let $\mu(f)=\int\mu(dx) f(x)$, be the Lebesgue
integral of a function $f\in\mathcal{B}(E)$, with respect to a measure
$\mu\in\mathcal{M}(E)$. We recall that a bounded integral operator $M$
from a measurable space $(E,\mathcal{E})$ into an auxiliary measurable
space $(F,\mathcal{F})$ is an operator $f\mapsto M(f)$ from
$\mathcal{B}(F)$ into $\mathcal{B}(E)$ such that the functions $
M(f)(x):=\int_{F}M(x,dy)f(y)$ are $\mathcal{E}$-measurable and bounded,
for any $f\in\mathcal{B}(F)$. A~Markov kernel is a positive and bounded
integral operator $M$ with $M(1)=1$. Given a pair of bounded integral
operators $(M_1,M_2)$, we let $(M_1M_2)$ the composition operator
defined by $(M_1M_2)(f)=M_1(M_2(f))$. For time homogenous state spaces,
we denote by $M^m=M^{m-1}M=MM^{m-1}$ the $m$th composition of a given
bounded integral operator $M$, with $m\geq1$. A~bounded integral
operator $M$ from a measurable space $(E,\mathcal {E})$ into an
auxiliary measurable space $(F,\mathcal{F})$ also generates a dual
operator $\mu\mapsto\mu M$ from $\mathcal{M}(E)$ into $\mathcal{M}(F)$
defined by $(\mu M)(f):=\mu(M(f))$. We also used the notation
\[
K\bigl([f-K(f)]^2\bigr)(x):=K\bigl( [f-K(f)(x) ]^2\bigr)(x)
\]
for some bounded integral operator $K$ and some bounded function $f$.

When the bounded integral operator $M$ has a constant mass, that is,
when $M(1)( x) =M(1)( y) $ for any $(x,y)\in E^{2}$, the operator
$\mu\mapsto\mu M$ maps $\mathcal{M}_{0}(E)$ into $\mathcal{M}_{0}(F)$.
In this situation, we let $\beta(M)$ be the Dobrushin coefficient of a
bounded integral operator $M$ defined by the formula
$\beta(M):=\sup\{\operatorname{osc}(M(f)) ; f\in
\operatorname{Osc}_{1}(F)\} $.

\subsection{Description of the main results}

The mathematical and numerical analysis of the mean field particle
models~(\ref{meanfield}) is one of the most attractive research areas
in pure and applied probability,
as well as in advanced
stochastic engineering and computational physics.

The fluctuation analysis of these discrete generation particle models
around their limiting distributions is often restricted to
Feynman--Kac-type models (see, e.g.,
\cite{chopine,fk,dmiclo,DG} and references therein) or specific
continuous time mean field models including McKean--Vlasov diffusions
and Boltzmann-type collision models of gases \cite{sylvie,Shiga}.

In the present article, we design an original and natural stochastic
perturbation analysis that applies to a rather large class of models
satisfying a rather weak first-order regularity property. We combine an
original stochastic perturbation analysis with a concentration analysis
for triangular arrays of conditionally independent random sequences,
which may be of independent interest. Under some additional stability
properties of the limiting measure valued processes, uniform
concentration properties with respect to the time parameter are also
derived. The concentration inequalities presented here generalize the
classical Hoeffding, Bernstein and Bennett inequalities for independent
random sequences to interacting particle systems, yielding very new
results for this class of models.

To describe with some precision this first main result
we observe that the local sampling errors associated with the
corresponding mean field particle model are expressed in terms of the
centered random fields $W_{n}^{N}$, given by
the following stochastic perturbation formulae:
%
%
\begin{equation}\label{defWNn}
\eta^N_n=\eta_{n-1}^NK_{n,\eta_{n-1}^N}+\frac{1}{\sqrt{N}} W^N_n
.
\end{equation}
To analyze the propagation properties of these local sampling errors,
\textit{up to a second-order remainder measure}, we further assume that
the one-step mappings $\Phi_n$ governing equation~(\ref{phi}) have a
first-order decomposition
%
%
\begin{equation}\label{decapp}
\Phi_{n}(\eta)-\Phi_{n}(\mu)\simeq(\eta-\mu) D_{\mu}{\Phi_n}
\end{equation}
with a first-order integral operator $D_{\mu}{\Phi_n}$ from $\Ba(E_n)$
into $\Ba(E_{n-1})$, s.t. $D_{\mu}{\Phi_n}(1)=0$.
The precise definition of the first-order regularity property~(\ref{decapp})
is provided in Definition~\ref{defregu}.

Our first main result is a functional central limit theorem for the
random fields
%
%
\begin{equation}\label{defVNn}
V^N_n:=\sqrt{N} [\eta^N_n-\eta_n] .
\end{equation}
This fluctuation theorem takes, basically, the following form.
\begin{theo}\label{tclinter}
\begin{itemize}
\item The sequence $(W^N_n)_{n\geq0}$
converges in law, as $N$ tends to infinity, to the sequence of $n$ independent,\vadjust{\goodbreak}
Gaussian and centered random fields $(W_n)_{n\geq0}$ with a covariance function
given for any
$f,g\in\Ba(E_n)$, and any $n\geq1$, by
%
%
\begin{equation}\label{corr1}\quad
\EE(W_n(f)W_n(g))=\eta_{n-1} K_{n,\eta_{n-1}}\bigl([f-K_{n,\eta_{n-1}}(f)][g-
K_{n,\eta_{n-1}}(g)]\bigr)
\end{equation}
and, for $n=0$, by
\[
\EE(W_0(f)W_0(g))=\eta_0\bigl[\bigl(f-\eta_0(f)\bigr)\bigl(g-\eta_0(g)\bigr)\bigr] .
\]
\item For any fixed time horizon $n\geq0$, the sequence of random fields
$V^N_n
$
converges in law, as the number of particles $N$ tends to infinity, to a
Gaussian and centered random fields
\[
V_n=\sum_{p=0}^nW_p \Da_{p,n} .
\]
In the above display, $\Da_{p,n}$ stands for the semigroup associated with
the operator $\Da_n=D_{\eta_{n-1}}{\Phi_n}$.
\end{itemize}
\end{theo}

A complete detailed proof of the functional central limit theorem
stated above
is provided in Section~\ref{stochpert}, dedicated to a
stochastic perturbation analysis of mean field particle
models.

We let $\Phi_{p,n}=\Phi_{p+1,n}\circ\Phi_{p+1} $, $0\leq p\leq n$, be
the semigroup associated with the measure valued equation defined in
(\ref{phi}). For $p=n$, we use the convention $\Phi_{n,n}=\mathrm{Id}$, the
identity operator.
By construction, we have
\[
\Phi_{p,n}(\eta)-\Phi_{p,n}(\mu)\simeq(\eta-\mu) D_{\mu}{\Phi _{p,n}}
\quad\mbox{and}\quad D_{\eta_p}{\Phi_{p,n}}=\Da_{p,n} .
\]
The fluctuation theorem stated above shows that the fluctuations of
$\eta^N_n$ around the limiting measure
$\eta_n$
is precisely dictated by
first-order differential-type operators
$D_{\eta_p}{\Phi_{p,n}}$ of the semigroup $\Phi_{p,n}$ around the flow
of measures $\eta_p$, with $p\leq n$.
Furthermore,
for any $f_n\in\operatorname{Osc}_1(E_n)$, one observes that
%
%
\begin{eqnarray}\label{asymptop}
\EE(V_n(f_n)^2)&=&\sum_{p=0}^n\EE((W_p[D_{\eta_p}
\Phi_{p,n}(f_n)])^2)\nonumber\\[-8pt]\\[-8pt]
&\leq&\sum_{p=0}^n\sigma_p^2
\beta(D\Phi_{p,n})^2:=\overline{\sigma}{}^2_n\nonumber
\end{eqnarray}
with the uniform local variance parameters
\[
\sigma_n^2:=\sup_{f_n\in\operatorname{Osc}_1(E_n)}{\sup_{\mu\in\Pa
(E_{n-1})}{\bigl| \mu\bigl( K_{n,\mu}[f_n-K_{n,\mu}(f_{n})]^2\bigr) \bigr|}}
\]
and
\[
\beta(D\Phi_{p,n}):=\sup_{\eta\in \Pa (E_p)}\beta(D_{\eta}\Phi_{p,n}) .
\]

The second part of this article is concerned with the concentration
properties of mean field particle
models.
These results quantify
exponentially small probabilities of deviations events between the
occupation measures $\eta^N_n$ and their\vadjust{\goodbreak} limiting values.
The exponential deviation events discussed in this article are
described in terms of the parameters
\[
\overline{\sigma}^2_n\leq\beta_n^2:= \sum_{p=0}^n\beta(D\Phi_{p,n})^2
\quad\mbox{and}\quad
b^{\star}_n:=\sup_{0\leq p\leq n}\beta(D\Phi_{p,n}) .
\]
Besides the fact that the
nonasymptotic analysis of weakly dependent variables is rather
well developed, the
concentration properties of discrete generation and interacting
particle systems often resume to
asymptotic large deviation results, or to nonasymptotic rough
exponential estimates (see, e.g., \cite{fk} and references
therein). Our main result on this subject is an original concentration
theorem that includes Hoeffding, Bennett and Bernstein
exponential inequalities for mean field particle models. This result
takes, basically, the following form.
\begin{theo}\label{conctheo}
For any $N\geq1$, $n\geq0$, $f_n\in\operatorname{Osc}_1(E_n)$, and any
$x\geq0$ the probability of each of the following pair of events
is greater than $1-e^{-x}$
\[
V^N_n(f_n)\leq\frac{r_n}{\sqrt{N}} \bigl(1+ \varepsilon_{0}^{-1} (x)\bigr)
+\sqrt{N} \overline{\sigma}_n^2 b_n^{\star}
\varepsilon_{1}^{-1}\biggl(\frac{x}{N\overline{\sigma}_n^2}\biggr)
\]
and
\[
V^N_n(f_n) \leq\frac{r_n}{\sqrt{N}} \bigl(1+ \varepsilon_{0}^{-1} (x)\bigr) +
\sqrt{2x} \beta_n.
\]
In the above display, $r_n$ stands for some parameter whose values only depend
on the amplitude of the second-order terms in the development
(\ref{decapp}), and the pair of functions $(\varepsilon_0,\varepsilon
_1)$ are
defined by
%
%
\begin{equation}\label{defepsi}\quad
\varepsilon_0(\lambda)=\tfrac{1}{2}\bigl(\lambda-\log{(1+\lambda
)}\bigr) ,\qquad
\varepsilon_1(\lambda) = (1+\lambda)\log{(1+\lambda)}-\lambda.
\end{equation}
Under additional stability properties of the semigroup associated with
the limiting model~(\ref{phidef}), the parameters $(\overline{\sigma
}_n,\beta_n ,b^{\star}_n,r_n)$ are uniformly bounded w.r.t.
the time parameter.
\end{theo}

A complete detailed proof of the functional concentration inequalities
stated in Theorem~\ref{conctheo}
is provided in Section~\ref{secmvp}. Some of the consequences of the
concentration inequalities stated above are provided in Section \ref
{concentsec}. To give a flavor of these results, using
a Bernstein-type concentration inequality we will check that
%
%
\begin{equation}\label{almosts}
\limsup_{N\rightarrow\infty} \log{\mathbb{P} \biggl(V^N_n(f_n)\geq
\frac{r_n}{\sqrt{N}}+\lambda\biggr)}\leq-\frac {\lambda
^2}{2(b_n^{\star} \overline{\sigma}_n)^2} .
\end{equation}
The detailed proof of this asymptotic estimate is provided on
page \pageref{corconcc}. This observation shows that this concentration
inequality is ``almost'' asymptotically sharp,
with a variance-type term whose values are pretty close to the exact limiting
variances presented in~(\ref{asymptop}). A~more precise asymptotic estimate
would require a refined moderate deviation analysis. We hope to discuss
these properties
in a forthcoming study.\vadjust{\goodbreak}

The outline of the rest of the article is as follows.
To motivate the present article, we have collected in Section \ref
{secillus} three
different classes of abstract mean field
particle models that can be studied using the fluctuation and the
concentration analysis developed in this article.

In Section~\ref{defFS}, we discuss the main regularity properties
used in our analysis. In Section~\ref{concentsec}, we illustrate the
impact of Theorem~\ref{conctheo}
with some more
Bennett and Hoeffding-type concentration properties, as well as
Bernstein-type concentration inequalities and uniform exponential
deviation properties w.r.t. the time parameter.
Section~\ref{stochpert} is mainly concerned with the detailed proofs
of the
theorems stated above.
We combine a natural stochastic perturbation analysis with nonlinear semigroup
techniques that allow us to describe both the fluctuations and the
concentration of the mean field measures in terms
of the local error random field models introduced in~(\ref{defWNn}).
The functional central limit theorem
is proved in Section~\ref{ftclf}. In Appendix~\ref{sconvex}, we provide
a preliminary convex
analysis including estimates of inverses of Legendre--Fenchel
transformations of
classical convex functions needed in this article.
In Section~\ref{lemmconc}, we prove a technical concentration lemma for
triangular arrays
of conditionally independent random variables. In Section~\ref{secmvp},
we apply this lemma to prove concentration
inequalities for mean field models.

\section{Some illustrative examples}\label{secillus}

\subsection{Feynman--Kac models}\label{secFKmod}

As mentioned in the \hyperref[intro]{Introduction}, the first prototype
model we have in
mind is a class of Feynman--Kac distribution flow equation arising
in a variety of application areas including stochastic engineering,
physics, biology and Bayesian statistics.
These models are defined in terms of a series
of bounded and positive integral operators $Q_n$ from $E_{n-1}$ into
$E_n$ with the following
dynamical equation:
%
%
\begin{equation}\label{fkflows}
\forall f_n\in\Ba(E_n)\qquad \eta_n(f_n)={\eta _{n-1}(Q_n(f_n))}/{\eta
_{n-1}(Q_{n}(1))}
\end{equation}
with a given initial distribution $\eta_0\in\Pa(E_0)$. To avoid
unnecessary technical discussions we simplify the analysis and we
assume that
\[
\forall n\geq0\qquad 0<\inf_{x\in E_{n}}G_n(x)\leq\sup_{x\in
E_{n}}G_n(x)<\infty\qquad \mbox{with } G_n(x):=Q_{n+1}(1)(x) .
\]
Rewritten in a slightly different way, we have
\[
\eta_n=\Phi_{n}(\eta_{n-1}):=\Psi_{n-1}(\eta_{n-1})M_n
\qquad\mbox{with } M_n(f_n)=Q_n(f_n)/Q_n(1)
\]
and the Boltzmann--Gibbs transformation $\Psi_{n}$ from $\Pa(E_{n}) $
into itself given by
\[
\forall f_n\in\Ba(E_{n})\qquad \Psi_{n}(\eta_{n})(f_n)={\eta
_{n}(G_{n}f_n)}/{\eta_{n}(G_{n})} .
\]
Using the ratio formulation~(\ref{fkflows}) of the semigroup, we will
check in Appendix~\ref{FeynKac} that
the first-order decomposition~(\ref{decapp}) is met\vadjust{\goodbreak}
with the first-order operator defined by
\[
D_{\mu}\Phi_{n}(f):=\frac{1}{\mu Q_n(1)} Q_{n}\bigl(f-\Phi_{n}(\mu
)(f)\bigr) .
\]

The nonlinear filtering model~(\ref{exref}), the particle absorption
model~(\ref{exrefab}) and
the Boltzmann--Gibbs distribution flow~(\ref{exrefabg}) can be
abstracted in this framework
by setting
\[
Q_{n+1}(x,dx')=G_n(x)M_{n+1}(x,dx').
\]

We leave the reader to check that this flow of measures satisfy the
recursive equation~(\ref{phi})
for any choice of Markov transitions given below:
%
%
\begin{equation}\label{fkmodK}\qquad
K_{n+1,\eta_{n}}(x,dy)= \varepsilon_{n}G_{n}(x)
M_n(x,dy)+\bigl(1-\varepsilon_{n}G_{n}(x)\bigr) \Phi _{n+1}(\eta _{n})(dy) .
\end{equation}
In the above displayed formula $\varepsilon_{n}$ stands for some
$[0,1]$-valued parameters that may depend on the current measure $\eta
_{n}$ and such that $\Vert\varepsilon_{n}G_{n}\Vert\leq1$. In this
situation, the mean field $N$-particle model associated with the
collection of Markov transitions~(\ref{fkmodK}) is a combination of
simple selection/mutation genetic transition
$\xi_n\leadsto\widehat{\xi}_n=(\widehat{\xi}^i_n)_{1\leq i\leq
N}\leadsto\xi_{n+1}$. During\vspace*{2pt} the selection stage, with probability
$\varepsilon_{n}G_{n}(\xi_{n}^{i})$, we set $\widehat{\xi
}_{n}^{i}=\xi_{n}^{i}$; otherwise, the particle\vspace*{1pt} jumps to a new
location, randomly drawn from the discrete distribution
$\Psi_{n}(\eta_{n}^{N})$. During the mutation stage, each of the
selected particles $\widehat {\xi} 
_{n}^{i}\leadsto\xi_{n+1}^{i}$ evolves according to the transition
$M_{n+1}$. If we set $\varepsilon_n=0$, the above particle model reduces
to the simple genetic-type model discussed in~(\ref{meanfield0}) in the
\hyperref[intro]{Introduction}.

\subsection{Gaussian mean field models}

The concentration analysis presented in this article is not restricted
to Feynman--Kac-type
models. It also applies to McKean-type models associated with a
collection of multivariate
Gaussian-type
Markov transitions on $E_n=\RR^d$, defined by
%
%
\begin{eqnarray}\label{transG}
K_{n,\eta}(x,dy)&=&\frac{1}{\sqrt{(2\pi)^d\operatorname{det}(Q_n)}} \nonumber\\[-8pt]\\[-8pt]
&&{}\times\exp{\biggl\{-\frac{1}{2}\bigl(y-d_n(x,\eta)\bigr)^{\prime} Q_n^{-1}
\bigl(y-d_n(x,\eta)\bigr)\biggr\}} \,dy ,\nonumber
\end{eqnarray}
with a nonsingular, positive and semi-definite covariance matrix $Q_n$
and some sufficiently regular drift mapping
$d_n\dvtx(x,\eta)\in\RR^d\times\Pa(\RR^d)\mapsto d(x,\eta)\in\RR^d$.
In Appendix~\ref{aGauss}, for $d=1$ we will check that any linear drift
function $d_n$ of the form $d_n(x,\eta)=a_n(x)+\eta(b_n) c_n(x)$, with
some measurable (and nonnecessarily bounded) function $a_n$, and some
pair of functions
$b_n$ and $c_n\in\Ba(\RR)$,
the first-order decomposition~(\ref{decapp}) is met
with the first-order operator defined by
\begin{eqnarray*}
D_{\mu}\Phi_{n}(f)(x)&:=&[K_{n,\mu}(f)(x)-\Phi_n(\mu)(f)]\\
&&{} +b_n(x)
\int\mu(dy) c_n(y) K_{n,\mu}(y,dz) f(z) \bigl(z-d_n(y,\mu)\bigr).
\end{eqnarray*}

In this context, the $N$-mean field particle model is given by the
following recursion:
\[
\forall1\leq i\leq N\qquad \xi^i_n=d_n(\xi^i_{n-1},\eta
^N_{n-1})+W^i_n ,
\]
where $(W^i_n)_{i\geq0 }$ is a collection of independent and
identically distributed
$d$-valued Gaussian random variables with covariance matrix $Q_n$. The
connection between these discrete generation models and the more traditional
continuous time McKean--Vlasov diffusion models is as follows. Consider
the partial differential equation
\[
\partial_t\mu_t=\frac{1}{2}
\sum_{i,j=1}^{d}\partial^2 _{x^{i},x^{j}}(\mu_t)-\sum
_{i=1}^{d}\partial
_{x^{i}}(b^{i}_t(\point, \eta_t)\mu_t) ,
\]
where $\mu_t$ is a probability measure on $\RR^d$, and $b_t$ some drift
term associated with
some interaction kernels $b'_t$ and given by
\[
b_t(x, \mu)=\int b'_t(x,x') \mu(dx') .
\]
Under appropriate regularity conditions, one can show that $\eta_t$ is the
marginal distribution at time $t$ of the law of the solution of the
nonlinear stochastic differential
equation
%
%
\begin{equation}\label{signalcont}
d\overline{X}_{t}=b_t(\overline{X}_t,\mu_t)\,dt+dB_t,
\end{equation}
where $B_{t}$ is a $d$-dimensional Brownian motion. These models have
been introduced in the late-1960s by McKean \cite{mckean}. The
convergence of the mean field particle model associated with the
diffusion~(\ref{signalcont}) has been deeply studied in the mid-1990s
by Bossy and Talay \cite{bossy1,bossy2}, M\'el\'eard \cite{sylvie} and
Sznitman \cite{Sznitman}. We also refer the reader to the more recent
treatments on McKean--Vlasov diffusion models by Bolley, Guillin and
Malrieu \cite {Bolley} and Bolley, Guillin and Villani \cite{Villani}.
Besides the fact that these continuous time probabilistic models are
directly connected to a rather large class of physical equations, to
get some computationally feasible solution, some kind of time
discretization scheme is needed. Mimicking traditional time
discretization techniques of deterministic dynamical systems, several
natural strategies can be used. For instance, we can use a Euler-type
discretization of the diffusion given by~(\ref{signalcont}) as follows:
\[
X^{\Delta}_{t_n}-X^{\Delta}_{t_{n-1}}=b_{t_{n-1}}(X^{\Delta
}_{t_{n-1}},\mu_{t_{n-1}}) \Delta+
(B_{t_n}-B_{t_{n-1}})
\]
on the time mesh $(t_n)_{n\geq0}$, with $(t_{n}-t_{n-1})=\Delta$, with
some initial random variable with
distribution
$\mu_0=\operatorname{Law}(X_0^{\Delta})$.
In this situation, the elementary transitions of the approximated
random states $X^{\Delta}_{t_n}$
are of the form~(\ref{transG}), with the identity covariance matrices
$Q_n=\mathrm{Id}$, and the drift functions
$d_n(x,\eta)=x+b_{t_{n-1}}(x,\eta)$. The refined
convergence analysis of these discrete
time approximation models for more general models, including granular
media equations
is developed Malrieu and Talay \cite{Malrieu2,Talay}.\vadjust{\goodbreak}

We mention that the semigroup derivation approach for functional
fluctuation theorems and concentration
inequalities developed in this article do not apply directly to any
nonlinear diffusion equations with general interaction kernels
$c_t$. The
semigroup derivation technique requires one to control recursively time
the integrability properties of the semigroup
associated with the first and second-order derivative terms. A~rather
crude sufficient condition is to assume that the drift
terms $d_n$ is of the form discussed in the discrete time model.

\subsection{A McKean model of gases}\label{gases}

We end this section with a mean field particle model arising in fluid
mechanics. We consider a
measurable state space $(S_n,\Sa_n)$ with a countably generated
$\sigma
$-field and an $(\Sa_n\otimes\Ea_n)$-measurable mapping $a_n$
be a from $(S_n\times E_n)$ into $\RR_+$ such that
$\int\nu_n(ds) a_n(s,x)=1$, for any $x\in E_n$, and some bounded
positive measure $\nu_n\in\Ma(S_n)$. To illustrate this model, we can
take a partition
of the state $E_n=\bigcup_{s\in S_n}A_s$ associated with a countable set
$S_n$ equipped with the counting measure $\nu_n(s)=1$ and set
$a_n(s,x)=1_{A_s}(x)$. We let $K_{n+1,\eta}$ be the McKean transition
defined by
%
%
\begin{equation}\label{nuvelos}
K_{n+1,\eta}(x,dy)=\int\nu_n(ds) \eta(du) a_n(s,u) M_{n+1}((s,x),dy).
\end{equation}
In the above displayed formula, $M_n$ stands for some Markov transition from
$(S_n\times E_n)$ into $E_{n+1}$. The discrete time version of McKean's
two-velocities model for
Maxwellian gases corresponds to the time homogenous model on $E_n=S_n=\{
-1,+1\}$
associated with the counting measure $\nu_n$ and the pair of parameters
\[
a_n(s,x)=1_{s}(x) \quad\mbox{and}\quad M_{n+1}((s,x),dy)=\delta
_{sx}(dy) .
\]
In this situation, the measure valued equation~(\ref{phi}) takes the
following quadratic form:
\[
\eta_{n+1}(+1)=\eta_{n}(+1)^2+\bigl(1-\eta_{n}(+1)\bigr)^2 .
\]
We leave the reader to write out the mean field particle interpretation
of this model.
For more details on this model, we refer to \cite{Shiga}. In
Appendix~\ref{McKeangas}, we will check that the first-order
decomposition~(\ref{decapp}) is met with the first-order operator
defined by
\begin{eqnarray*}
D_{\mu}\Phi_{n+1}(f)(x)&=&[K_{n+1,\mu}(f)(x)-\Phi_{n+1}(\mu
)(f)]\\
&&{}+\int\nu_n(ds) [a(s,x)-\mu(a(s,\point))]
\mu(M_{n+1}(f)(s,\point)).
\end{eqnarray*}

\section{Some weak regularity properties}\label{defFS}

To describe precisely the concentration inequalities developed in the
article, we need to
introduce a first round of notation.\vadjust{\goodbreak}
\begin{defi}\label{defregu}
We let $\Upsilon(E,F)$ be the set of mappings $\Phi\dvtx\mu\in\Pa
(E)\mapsto\Phi(\mu)\in\Pa(F)$ satisfying the first-order decomposition
%
%
\begin{equation}\label{firstodef}
\Phi(\mu)-\Phi(\eta)=(\mu-\eta) D_{\eta}\Phi+\Ra^{\Phi}(\mu
,\eta) .
\end{equation}
In the above displayed formula, the first-order operators $(\Da_{\eta
}\Phi)_{\eta\in\Pa(E)}$ is some collection of bounded
integral operators from $E$ into $F$ such that
%
%
\begin{eqnarray}\label{firsto}
&\displaystyle \forall\eta\in\Pa(E), \forall x\in E\qquad
(D_{\eta}\Phi)(1)(x)=0 \quad\mbox{and}&\nonumber\\[-8pt]\\[-8pt]
&\displaystyle \beta(\Da\Phi)
:=\sup_{\eta\in\Pa(E)}\beta(D_{\eta}\Phi)<\infty.&\nonumber
\end{eqnarray}
The collection of second-order remainder signed measures $(\Ra^{\Phi
}(\mu,\eta))_{(\mu,\eta)\in\Pa(E^2)}$ on $F$ are such that
%
%
\begin{equation}\label{condxi222}
|\Ra^{\Phi}(\mu,\eta)(f)|\leq\int|(\mu-\eta
)^{\otimes
2}(g)| R^{\Phi}_{\eta}(f,dg) ,
\end{equation}
for some collection of integral operators $R^{\Phi}_{\eta}$ from $\Ba
(F)$ into the set $\operatorname{Osc}_{1}(E)^2$ such that
%
%
\begin{eqnarray}\label{lipr}
\sup_{\eta\in\Pa(E)} \int\operatorname{osc}(g_1)
\operatorname{osc}(g_2) R^{\Phi}_{\eta}\bigl(f,d(g_1\otimes g_2)\bigr)\leq
\operatorname{osc}(f) \delta(R^{\Phi})\nonumber\\[-8pt]\\[-8pt]
&&\eqntext{\mbox{with } \delta(R^{\Phi})<\infty.}
\end{eqnarray}
\end{defi}

This rather weak first-order regularity property is satisfied for a
large class of
one-step transformations $\Phi_n$ associated with a nonlinear measure valued
process~(\ref{phi}).
For instance, in Section~\ref{secmvp} we shall prove that the
Feynman--Kac transformations $\Phi_n$
introduced in~(\ref{fkflows}) belong to the set $ \Upsilon
(E_{n-1},E_n)$. The latter is also met for the
Gaussian transitions introduced in~(\ref{transG}) and for the
McKean-type model of gases~(\ref{nuvelos})
presented in Section~\ref{gases}.
The proof of this assertion is rather technical and it is postponed in
Appendix~\ref{McKeangas}.

We assume that the one-step mappings
\[
\Phi_n \dvtx\mu\in\Pa(E_{n-1})\longrightarrow\Phi_n(\mu):=\mu
K_{n,\mu
}\in\Pa(E_n)
\]
governing equation~(\ref{phi})
are chosen so that $\Phi_n\in\Upsilon(E_{n-1},E_n)$, for any $n\geq1$.
The main advantage of the regularity condition comes from the fact that
$\Phi_{p,n}\in\Upsilon(E_{p},E_n)$ with
the first-order decomposition-type formula
\[
\Phi_{p,n}(\eta)-\Phi_{p,n}(\mu)=[\eta-\mu]D_{\mu}\Phi
_{p,n}+\Ra^{\Phi
_{p,n}}(\eta,\mu) ,
\]
for some collection of bounded integral operators $D_{\mu}\Phi_{p,n}$
from $E_p$ into $E_n$ and
some second-order remainder signed measures $\Ra^{\Phi_{p,n}}(\eta
,\mu
)$. For further use, we let
$r_n$ be the second-order stochastic perturbation term related to the
quadratic remainder measures $R^{\Phi_{p,n}}$ and defined by
\[
r_n:= \sum_{p=0}^n\delta(R^{\Phi_{p,n}}) .
\]

\section{Some exponential concentration inequalities}\label{concentsec}

Let us examine some more or less direct consequences of the
concentration inequalities
stated in Theorem~\ref{conctheo}.

When the Markov kernels $K_{n,\mu}=K_n$ do not depend on the measure
$\mu$, the $N$-particle model reduces to a collection of independent
copies of the Markov chain with elementary transitions $P_n=K_n$.
In this special case, the second-order parameters vanish (i.e.,
$r_n=0$), while
the first-order expansion parameters $(\overline{\sigma}_n, \beta_n)$
are related to
the mixing properties of the semigroup
of the underlying Markov chain; that is, we have that
\[
\overline{\sigma}^2_n=\sum_{p=0}^n\sigma_p^2 \beta(P_{p,n})^2\leq
\beta_n^2= \sum_{p=0}^n\beta(P_{p,n})^2 \qquad\mbox{with }
P_{p,n}=K_{p+1},\ldots, K_{n-1}K_n ,
\]
with the Dobrushin ergodic coefficient $\beta(P_{p,n})$ associated with
$P_{p,n}$.
When the chain is asymptotically stable in the sense that $\sup_{n\geq
0} \sum_{p=0}^n\beta(P_{p,n})<\infty$,
the first-order expansion parameters given above are uniformly bounded
with respect to
the time parameter.

In more general situations, the analysis of these parameters depends
on the model at hand. For instance, for time homogeneous Feynman--Kac
models [i.e., $E_n=E$ and $(G_n,M_n)=(G,M)$]
these parameters can be related to the mixing properties of the Markov
chain associated with the transitions~$M$.
To be more precise, let us suppose that the following condition is
met:\looseness=-1
%
\begin{eqnarray}\label{fkhom1}\hypertarget{fkhom1link}{}
&(M)_m \quad\exists m\geq1, \exists\varepsilon_m>0
\quad\mbox{s.t.}&\nonumber\\[-8pt]\\[-8pt]
&\displaystyle \forall(x,y)\in E^{2} \qquad M^m%
(x,\cdot)\geq\varepsilon_m M^m(y,\cdot) .&\nonumber
\end{eqnarray}\looseness=0
It is well known that the mixing-type condition \hyperlink{fkhom1link}{$(M)_m$} is
satisfied for any aperiodic and irreducible Markov chains on finite
spaces, as
well as for bi-Laplace exponential transitions associated with a
bounded drift
function and for Gaussian transitions with a mean drift function that
is constant outside some
compact domain. To go one step further, we introduce the following quantities:
%
%
\begin{equation}\label{fkhom2}
\delta_{m}:=\sup\prod_{0\leq p<m} \bigl({G(x_{p})}/{G(y_{p})}\bigr).
\end{equation}
In the above displayed formula, the supremum is taken over all
admissible pair of paths with elementary transitions $M$. In this
situation, we can check that
\[
r_n\leq4 \varpi_{3,1}(m),\qquad
b^{\star}_n\leq2 \delta_{m}/\varepsilon_m
\]
as well as
\[
\overline{\sigma}^2_n\leq4 \varpi_{2,2}(m) \sigma^2
\quad\mbox{and}\quad
\beta_n^2\leq4 \varpi_{2,2}(m)
\]
with the uniform local variance parameter $\sigma^2$ and a collection
of parameters
$\varpi_{k,l}(m)$ such that $\varpi_{k,l}(m)\leq m \delta
_{m-1} \delta
_{m}^k/\varepsilon_m^{k+2}$.
The detailed proof of these estimates can be found in
Appendix~\ref{FeynKac}.\vadjust{\goodbreak}

As we mentioned above, in the special case where the Markov kernels
$K_{n,\mu}=K_n$ do not depend on the measure
$\mu$, the random measures $\eta^N_n$ coincide with the occupation
measure associated with $N$ independent and identically distributed
random variables with common law $\eta_n$. In this situation, the pair
of events described in Theorem~\ref{conctheo} resumes
to the following Bennett and Hoeffding-type concentration events,
respectively, given by
\[
[\eta^N_n-\eta_n](f_n)\leq\overline{\sigma}_n^2 b_n^{\star}
\varepsilon_{1}^{-1}\biggl(\frac{x}{N\overline{\sigma}_n^2}\biggr)
\quad\mbox{and}\quad
[\eta^N_n-\eta_n](f_n)\leq\sqrt{\frac{2x}{N}} \beta_n .
\]
The first inequality can be described more explicitly using the
analytic estimates
\[
\varepsilon_1^{-1}(x)\leq\frac{ \sqrt{2x} + (4x/3) - \log(1 +(x/3) +
\sqrt{2x})}{\log(1 +(x/3) + \sqrt{2x})} \leq(x/3) + \sqrt{2x} .
\]

In the context of Feynman--Kac models, the second-order terms can be
estimated more explicitly
using
the upper bounds
\[
\varepsilon_0^{-1} (x)\leq2x+\log\bigl( 1+2x+ 2\sqrt{x} \bigr) +
\frac{ \log( 1+2x+ 2\sqrt{x} ) - 2 \sqrt{x} }{2 x+ 2 \sqrt{x} }
\leq2 x+ 2 \sqrt{x} .
\]
A detailed proof of the upper bounds given above is detailed in
Appendix~\ref{sconvex}, dedicated to the convex
analysis of the Legendre--Fenchel transformations used in this article.
The second rough estimate in the r.h.s.
of the above displayed formulae leads to Bernstein-type concentration
inequalities.
\begin{cor}\label{corconcc}
For any $N\geq1$ and any $n\geq0$, we have the following
Bernstein-type concentration inequalities:
\begin{eqnarray*}
\frac{1}{N} \log{\mathbb{P}\biggl([\eta^N_n-\eta_n](f_n)\geq\frac
{r_n}{N}+\lambda\biggr)} &\geq& \frac{\lambda^2}{2} \biggl(\biggl(
b_n^{\star} \overline{\sigma}_n+\frac{\sqrt{2}r_n}{\sqrt{N}}
\biggr)^2+\lambda\biggl(2r_n+\frac{b_n^{\star}}{3}\biggr)\biggr)^{-1}
\end{eqnarray*}
and
\begin{eqnarray*}
-\frac{1}{N} \log{\mathbb{P}\biggl([\eta^N_n-\eta_n](f_n)\geq\frac
{r_n}{N}+\lambda\biggr)}&\geq& \frac{\lambda^2}{2} \biggl(
\biggl(\beta_n+ \frac{\sqrt{2}r_n}{\sqrt{N}} \biggr)^2
+2r_n\lambda\biggr)^{-1} .
\end{eqnarray*}
\end{cor}

In terms of the random fields $V^N_n$, the first concentration inequality
stated in Corollary~\ref{corconcc} takes the following form:
\begin{eqnarray*}
&&- \log{\mathbb{P}\biggl(V^N_n(f_n)\geq\frac{r_n}{\sqrt
{N}}+\lambda
\biggr)}\\
&&\qquad\geq\frac{\lambda^2}{2} \biggl(\biggl(
b_n^{\star}
\overline{\sigma}_n+\frac{\sqrt{2}r_n}{\sqrt{N}}
\biggr)^2+\frac{\lambda}{\sqrt{N}} \biggl(2r_n+\frac{b_n^{\star
}}{3}\biggr)\biggr)^{-1}\\
&&\qquad\mathop{\longrightarrow}\limits_{N\rightarrow\infty}
\frac{\lambda^2}{2(b_n^{\star}
\overline{\sigma}_n)^2} .
\end{eqnarray*}
This proves the asymptotic estimate presented in~(\ref{almosts}).\vadjust{\goodbreak}

Last, but not least, without further work, Theorem~\ref{conctheo} leads
to uniform concentration inequalities for
mean field particle interpretations of Feynman--Kac semigroups.
\begin{cor}
In the context of Feynman--Kac models, under the mixing type
condition \hyperlink{fkhom1link}{$(M)_m$} introduced in~(\ref{fkhom1}), for any $N\geq1$, any
$n\geq0$ and any $x\geq0$ the probability of each of the following
pair of events:
\begin{eqnarray*}
[\eta^N_n-\eta_n](f_n)
&\leq&\frac{4}{N}\varpi_{3,1}(m)\bigl(1+
\varepsilon_{0}^{-1} (x)\bigr)\\
&&{}+\frac{8\delta_{m}}{\varepsilon_m} \varpi_{2,2}(m) \sigma^2
\varepsilon_{1}^{-1}\biggl(\frac{x}{4\sigma^2\varpi_{2,2}(m) N}\biggr)
\end{eqnarray*}
and
\[
[\eta^N_n-\eta_n](f_n)\leq\frac{4}{N}\varpi_{3,1}(m)\bigl(1+
\varepsilon_{0}^{-1} (x)\bigr) +2 \sqrt{\frac{2\varpi_{2,2}(m)x}{N}}
\]
is greater than $1-e^{-x}$.
\end{cor}

\section{A stochastic perturbation analysis}\label{stochpert}
\subsection{Proof of the functional central limit theorem}\label{ftclf}
\begin{defi}
We say that a collection of
Markov transitions $K_{\eta}$ from a measurable space $(E,\Ea)$ into
another $(F,\Fa)$
satisfies condition \hyperlink{lipsKlink}{$(K)$} as soon as the
following Lipschitz-type inequality is met for every $f\in\operatorname
{Osc}_1(F)$:
%
%
\begin{equation}\label{lipsK}\hypertarget{lipsKlink}{}
(K)\quad\Vert[ K_{\mu}-K_{\eta}](f)\Vert\leq
\int|(\mu
-\eta)(h)| T^K_{\eta}(f,dh) .
\end{equation}
In the above display, $T^K_{\eta}$ stands for some collection of
bounded integral operators from $\Ba(F)$ into $\Ba(E)$ such that
%
%
\begin{equation}
\sup_{\eta\in\Pa(E)} \int\operatorname{osc}(h) T^K_{\eta
}(f,dh)\leq
\operatorname{osc}(f) \delta(T^{K}) ,
\end{equation}
for some finite constant $\delta(T^{\Phi})<\infty$. In the
special case where $K_{\eta}(x,dy)=\Phi(\eta)(dy)$, for some mapping
$\Phi\dvtx\eta\in\Pa(E)\mapsto\Phi(\eta)\in\Pa(F)$, condition
(\ref{lipsK})
is a simple Lipschitz-type condition on the mapping $\Phi$. In this
situation, we denote by
\hyperlink{PhiLink}{$(\Phi)$} the corresponding condition; and whenever it is met, we says
that the mapping $\Phi$ satisfy condition \hyperlink{PhiLink}{$(\Phi)$}.
\end{defi}

We further assume that we are given a collection of McKean transitions
$K_{n,\eta}$ satisfying the weak Lipschitz-type condition stated in
(\ref{lipsK}). In this situation, we already mention that
the corresponding one-step mappings $\Phi_n(\eta)=\eta K_{n,\eta}$, and
the corresponding
semigroup $\Phi_{p,n}$ satisfies condition $(\Phi_{p,n})$ for some
collection of bounded integral operators $T^{\Phi_{p,n}}_{\eta}$.\vadjust{\goodbreak}

In the context of Feynman--Kac-type mdels, it is not difficult to check
that condition $(\Phi_n)$ is equivalent to
the fact that the McKean transitions $K_{n,\eta}$ given in (\ref
{fkmodK}) satisfy the Lipschitz condition~(\ref{lipsK}). The latter is
also met for the Gaussian transitions introduced in~(\ref{transG}) as
soon as
the drift function $d(x,\eta)$ is sufficiently regular. As before, this
condition is met for the Gaussian transitions introduced in (\ref
{transG}) and for the McKean-type model of gases~(\ref{nuvelos})
presented in Section~\ref{gases}.
For a more detailed discussion on these stability properties, we refer
the reader to the \hyperref[app]{Appendix}, on page 23.

Notice that the centered random fields $W^N_n$ introduced in (\ref
{defWNn}) have
conditional variance functions given by
%
%
\begin{equation}\label{covf}
\EE(W_{n}^{N}(f_{n})^2| \Fa^{N}_{n-1} )=\eta
_{n-1}^N
\bigl[K_{n,\eta
_{n-1}^{N}}\bigl(\bigl(f_n-K_{n,\eta
_{n-1}^{N}}(f_{n})\bigr)^2\bigr)\bigr] .
\end{equation}
Using Kintchine's inequality,
for every $f\in\operatorname{Osc}_1(E_n)$, $N\geq1$ and any $n\geq0$
and $m\geq1$ we have
the $\LL_{2m}$ almost sure estimates
%
%
\begin{equation}\label{kint}
\EE\bigl(
|W_{n}^{N}(f_{n})|^{2m}
| \Fa^{(N)}_{n-1}
\bigr)^{{1}/({2m})}\leq b(2m) \quad\mbox{with }
b(2m)^{2m}:=2^{-m} (2m)!/m!.\hspace*{-35pt}
\end{equation}

We can also prove the following theorem.
\begin{theo}\label{theocentralfluc}
The sequence $(W^N_n)_{n\geq0}$
converges in law, as $N$ tends to infinity, to the sequence of $n$ independent,
Gaussian and centered random fields $(W_n)_{n\geq0}$
described in Theorem~\ref{tclinter}.
\end{theo}

The proof of this theorem follows the same line of arguments as those
we used in \cite{fk} in the context of Feynman--Kac models. For
completeness, and for the convenience of the reader, the complete proof
of this result is housed in Appendix~\ref{secappfluc}.

Let us examine some direct consequences of this result.
Combining the Lipschitz property $(\Phi_{p,n})$ of the semigroup $\Phi
_{p,n}$ with the decomposition
\[
[\eta^N_n-\eta_n]=\sum_{p=0}^n [ \Phi_{p,n}(\eta^N_p)-
\Phi_{p,n}(\Phi_{p}(\eta^N_{p-1})) ] ,
\]
we find that
\[
\sqrt{N} |[\eta^N_n-\eta_n](f_n)|=\sum_{p=0}^n
\int|W^N_p(h)| T^{\Phi_{p,n}}_{\Phi_{p}(\eta
^N_{p-1})}(f,dh) .
\]
In the above displayed formulae, we have used the convention $\Phi
_{0}(\eta^N_{-1})=\eta_0$, for $p=0$. From the previous $\LL_{2m}$
almost sure estimates, we readily conclude that
\[
\sup_{N\geq1}{\sqrt{N} \EE\bigl(|[\eta^N_n-\eta_n
](f_n)|^{2m} \bigr)^{{1}/({2m})}}\leq b(2m) \sum_{p=0}^n\delta
(T^{\Phi_{p,n}}) .
\]

We are now in position to prove the fluctuation Theorem~\ref{tclinter}.
Using the decomposition
\[
V^N_n=W^N_n+V^N_{n-1}\Da_n+\sqrt{N} R^{\Phi_n}(
\eta^N_{n-1},\eta_{n-1} ) ,
\]
we readily prove that
%
%
\begin{equation}\label{decompV}
V^N_n=\sum_{p=0}^nW^N_p\Da_{p,n}+\frac{1}{\sqrt{N}} \Ra^N_n ,
\end{equation}
with the remainder second-order measure
\[
\Ra^N_n:=
N \sum_{p=0}^{n-1}R_{p+1}^{\Phi_{p+1}}(\eta^N_{p},\eta
_p
)D_{p+1,n} .
\]
In the above display, $\Da_{p,n}=\Da_{p+1},\ldots,\Da_{n-1}\Da_n$ stands
for the semigroup associated with the integral operators $\Da
_n:=D_{\eta
_{n-1}}\Phi_n$,
with the usual convention $\Da_{n,n}=\mathrm{Id}$, for $p=n$. Using a
first-order derivation formula for
the semigroup $\Phi_{p,n}$ (cf., e.g., Lemma~\ref{derivform}
on page \pageref{derivform}),
it is readily checked that
\[
D_{\eta_p}\Phi_{p,n}=(D_{\eta_p}\Phi_{p+1})(D_{\eta_{p+1}}\Phi
_{p+1,n})=\Da_{p+1}(D_{\eta_p}\Phi_{p,n})=\Da_{p,n} .
\]
Using the fact that
\[
|\Ra^N_n(f_n)|\leq
\sum_{p=0}^{n-1}
\int|(V^N_p)^{\otimes2}(g)| R^{\Phi
_{p+1}}_{\eta
_{p}}(f,dg) ,
\]
we conclude that, for any $m\geq1$, we have
\[
\EE(|\Ra^N_n(f_n)|^m)^{1/m}\leq
b(2m)^2 \sum_{p=0}^{n-1}\beta(\Da_{p+1,n}) \Biggl(\sum
_{q=0}^p\delta
(T^{\Phi_{q,p}})\Biggr)^2 \delta(R^{\Phi_{p+1}}) .
\]
This clearly implies that $\frac{1}{\sqrt{N}} \Ra^N_n$ converge in law
to the null measure, in the sense
that $\frac{1}{\sqrt{N}} \Ra^N_n(f_n)$ converge in law to zero, for any
bounded test function $f_n$ on $E_n$. Using the fact that $W^N_n$
converges in law to the sequence of $n$ independent,
random fields $W_n$, the proposition is now a direct consequence of the
decomposition formula~(\ref{decompV}). This ends the proof of
Theorem~\ref{tclinter}.

\subsection{A concentration lemma for triangular arrays}\label{lemmconc}

For every $n\geq0$ and $N\geq1$, we let
$X^{(N)}_n:=(X^{(N,i)}_n)_{1\leq i\leq N}$ be a triangular array of
random variables defined on some
filtered probability space $(\Omega,\Fa^N_n)$ associated with a
collection of increasing $\sigma$-fields $(\Fa^N_n)_{n\geq0}$. We
assume that $(X^{(N,i)}_n)_{1\leq i\leq N}$ are $\Fa
^N_{n-1}$-conditionally independent and centered random variables.
Suppose furthermore that
\[
\forall n\geq0\qquad
a_n\leq X^{(N,i)}_n\leq b_n
\quad\mbox{and}\quad
\EE\bigl(\bigl(X^{(N,i)}_n\bigr)^2 | \Fa^N_{n-1}\bigr)\leq c_n^2
\]
for some collection of finite constants $(a_n,b_n,c_n)$, with the
convention $\Fa^N_{-1}=\{\varnothing,\Omega\}$ for $n=0$. For any
$n\geq
0$, let
\[
T^N_n:=S^N_n+R_n^N \qquad
\mbox{where } \Delta S_n^N:=S_{n}^N-S^N_{n-1}=\sum_{i=1}^N
X^{(N,i)}_n
\]
and $R_n^N$ is a random perturbation term such that
\[
\forall m\geq1\qquad \EE(|R_n^N|^m)^{
{1/m}}\leq b(2m)^{2} d_n
\]
for some finite constant $d_n$. We use the convention $S^N_{-1}=0$, for
$n=0$. We set
\[
\overline{c}_n^2:=(b^{\star}_n)^{-2}\sum_{p=0}^nc_p^2
\quad\mbox{and}\quad
\overline{\delta}{}^2_n:=\sum_{p=0}^n\delta_p^2
\]
with the middle point
\[
\delta_n:=\frac{b_n-a_n}{2} .
\]

\begin{lem}\label{lemm}
For any $N\geq1$ and any $n\geq0$, the probability of each of the
following pair of events
%
%
\begin{equation}\label{eqalpha}
T^N_n\leq d_n\bigl(1+
\varepsilon_{0}^{-1} (x)\bigr)
+N \overline{c}_n^2 b_n^{\star}
\varepsilon_{1}^{-1}\biggl(\frac{x}{N\overline{c}_n^2}\biggr)
\end{equation}
and
%
%
\begin{equation}\label{eqalpha2}
T^N_n\leq d_n\bigl(1+
\varepsilon_{0}^{-1} (x)\bigr)
+\overline{\delta}_n
\sqrt{2xN}
\end{equation}
is greater than $1-e^{-x}$, for any $x\geq0$.
\end{lem}
\begin{rem} Notice that~(\ref{eqalpha2}) gives always a better
concentration inequality when
$\sum_{p=0}^nc_p^2\geq\sum_{p=0}^n\delta_p^2$.
In the opposite situation, if $\sum_{p=0}^nc_p^2 < \sum_{p=0}^n\delta_p^2$,
inequality~(\ref{eqalpha}) gives better concentration estimates for
sufficiently
small values of the precision parameter $x$.
\end{rem}

Before getting into the details of the proof of the above lemma, we
examine some
direct consequences of these inequalities based on Legendre--Fenchel
transforms estimates
developed in Appendix~\ref{sconvex}.
First, combining~(\ref{majstar0}) with~(\ref{majstar1}) we observe that,
with probability greater than $1-e^{-x}$,
\[
T^N_n\leq d_n\bigl(1+2\sqrt{x} +\theta_0(x)\bigr) +b_n^{\star}
\biggl(\overline{c}_n \sqrt{N} \sqrt{2x} + N\overline{c}_n^2
\theta_{1}\biggl(\frac{x}{N\overline{c}_n^2}\biggr) \biggr)
\]
with the pair of functions
\[
\theta_0(x):= 2x + \log\bigl(1+2\sqrt{x}+2x\bigr) - 2\sqrt{x} +
\frac{ \log(1+2\sqrt{x} +2x)) - 2\sqrt{x} }{2x + 2\sqrt{x} }
\leq2x
\]
and
\[
\theta_1(x):=\frac{\sqrt{2x} + (4x/3)}{\log(1+(x/3) + \sqrt{2x})}
- 1
- \sqrt{2x}
\leq\frac{x}{3} .
\]
The upper bounds given above together with~(\ref{breta}) imply that,
with probability greater than $1-e^{-x}$,
\[
T^N_n \leq d_n+A_n x +\sqrt{2x B_n^N} ,
\]
where
\[
A_n:=\biggl(2d_n+\frac{b_n^{\star}}{3}\biggr)
\quad\mbox{and}\quad
B_n^N:=\bigl( \sqrt{2}d_n+b_n^{\star} \overline{c}_n \sqrt{N}
\bigr)^2 .
\]
Using these successive upper bounds, we arrive at the following
Bernstein-type inequality:
%
%
\begin{eqnarray}\label{ineq1}
&&-\frac{1}{N} \log\mathbb{P}\biggl(\frac{T^N_n}{N}\geq\frac
{d_n}{N}+\lambda
\biggr)\nonumber\\[-8pt]\\[-8pt]
&&\qquad\geq
\frac{\lambda^2}{2} \biggl(\biggl(b_n^{\star}
\overline{c}_n+\frac{\sqrt{2}d_n}{\sqrt{N}}
\biggr)^2+\lambda\biggl(2d_n+\frac{b_n^{\star}}{3}\biggr)
\biggr)^{-1} .\nonumber
\end{eqnarray}

In much the same way, starting from~(\ref{eqalpha}), we have, with
probability greater
than $1-e^{-x}$,
%
%
\begin{equation}
T^N_n \leq d_n\bigl(1+ 2\bigl(x+\sqrt{x}\bigr)\bigr)
+\overline{\delta}_n
\sqrt{2xN}=d_n+A_n x
+\sqrt{2x B_n^N} ,
\end{equation}
with the pair of constants
\[
A_n:=2d_n \quad\mbox{and}\quad
B_n^N:=\bigl( \sqrt{2}d_n+\overline{\delta}_n \sqrt{N} \bigr)^2 .
\]
Using these successive upper bounds, we arrive at the following
Bernstein-type inequality:
%
%
\begin{equation} \label{ineq2}\qquad
-\frac{1}{N} \log{\mathbb{P}\biggl(\frac{T^N_n}{N}\geq\frac
{d_n}{N}+\lambda
\biggr)} \geq
\frac{\lambda^2}{2} \biggl( \biggl(\overline{\delta}_n+ \frac
{\sqrt
{2}d_n}{\sqrt{N}}
\biggr)^2 +2d_n\lambda\biggr)^{-1} .
\end{equation}
\begin{pf*}{Proof of Lemma~\ref{lemm}}
First, we observe that
\[
\forall t\in[0,1/(2d_n)[\qquad \EE
(e^{tR^N_n})\leq
\sum_{m\geq0}\frac{(td_n)^m}{m!} b(2m)^{2m} .
\]
To obtain a more explicit form of the r.h.s. term, we recall that
$b(2m)^{2m}=\EE(X^{2m})$ with a Gaussian centered random variable with
$\EE(X^2)=1$ and
\[
\forall d\in[0,1/2[\qquad \EE(\exp{\{dX^2\}})=\sum
_{m\geq
0}\frac{s^m}{m!}
b(2m)^{2m}= \frac{1}{\sqrt{1-2d}} .
\]
From this observation, we readily find that
\[
\forall t\in[0,1/(2d_n)[\qquad L^{N}_{0,n}(t):=\log{\EE
\bigl(e^{t(R^N_n-d_n)}\bigr)}\leq\alpha_{0,n}(t):=\alpha_0(td_n) .\vadjust{\goodbreak}
\]

Using~(\ref{ben}), we obtain the following almost sure inequality:
\[
\log{\EE(e^{t\Delta S^N_n} | \Fa^N_{n-1} )}\leq N
\biggl(\frac{c_n}{b_n}\biggr)^2 \alpha_1( b_nt ) .
\]
It implies that
\[
\forall t\geq0\qquad L^{N}_{1,n}(t):=\log{\EE(e^{tS^N_n}
)}\leq N \sum_{p=0}^n\biggl(\frac{c_p}{b_p}\biggr)^2 \alpha_1(
b_pt )\leq\alpha_{1,n}^N(t) ,
\]
with the increasing and convex function
$\alpha_{1,n}^N(t)=N \overline{c}_n^2 \alpha_1(
b_n^{\star}t)$.

Using~(\ref{breta}), we now obtain the following Cram\'{e}r--Chernoff estimate:
%
%
\begin{equation}\label{eqrio}\qquad
\forall x\geq0\qquad
\mathbb{P}\bigl(S^N_n+R^N_n\geq r_n+(L^{N \star}_{0,n}
)^{-1}(x)+(L^{N \star}_{1,n})^{-1}(x)\bigr)
\leq e^{-x} .
\end{equation}
In other words, the probability that
\[
S^N_n+R^N_n\leq r_n+(L^{N \star}_{0,n})^{-1}(x)+
(L^{N \star}_{1,n})^{-1}(x)
\]
is greater than $1-e^{-x}$, which, together with the homogeneity
properties of
the inverses of Legendre--Fenchel transforms recalled in Appendix~\ref{sconvex},
gives~(\ref{eqalpha}).

The proof of~(\ref{eqalpha2}) is based on Hoeffding's inequality,
\[
8 \log{\EE(e^{tX^{(N,i)}_n} | \Fa_{n-1}^N
)}\leq
t^2 (b_n-a_n)^2.
\]
From these estimates, we readily find that
$L^{N}_{1,n}(t)\leq\alpha_{2,n}^N(t):=N \overline{\delta}{}^2_n t^2/2$.
Arguing as before, we find that
\[
(L^{N \star}_{1,n})^{-1}(x)\leq(\alpha^{N \star
}_{2,n})^{-1}(x) =\sqrt{2xN\overline{\delta}{}^2_n}.
\]
We end the proof of the second assertion using~(\ref{eqrio}).
This ends the proof of the lemma.
\end{pf*}

\subsection{Concentration properties of mean field models}\label{secmvp}

This section is concerned with the proof of Theorem~\ref{conctheo}. To
simplify the presentation, we set
\[
\Da_{p,n}^{(N)}:=\Da_{\Phi_{p}(\eta^N_{p-1}) }\Phi_{p,n}
\quad\mbox{and}\quad
\Ra_{p,n}=\Ra^{\Phi_{p,n}} .
\]
Under our assumptions, we have the almost sure estimates
\[
\sup_{N\geq1}\beta\bigl(\Da_{p,n}^{(N)}\bigr)\leq\beta
(\Da\Phi
_{p,n}):=\sup_{\eta\in\Pa(E_p)}\beta(\Da_{\eta}\Phi
_{p,n}) .
\]
In this notation, one important consequence of the above lemma is the
following decomposition:
\begin{eqnarray*}
V^N_n:\!&=&\sqrt{N} [\eta^N_n-\eta_n]\\
&=&\sqrt{N} \sum_{p=0}^n [ \Phi_{p,n}(\eta^N_p)- \Phi
_{p,n}
(\Phi_{p}(\eta^N_{p-1})) ]=I^N_n+J^N_n
\end{eqnarray*}
with the pair of random measures $(I^N_n,J^N_n)$ given by
\begin{eqnarray*}
I^N_n:=\sum_{p=0}^nW^N_p\Da_{p,n}^{(N)} \quad\mbox{and}\quad
J^N_n:=\sqrt{N} \sum_{p=0}^n\Ra_{p,n}(\eta^N_p,\Phi_{p}(\eta
^N_{p-1})) .
\end{eqnarray*}
In what follows $f_n$ stands for some test function $f_n\in\mbox
{Osc}_{1}(E_n)$.
Combining~(\ref{kint}) with the generalized Minkowski integral
inequality we find that
\[
N \EE\bigl(|\Ra_{p,n}(\eta^N_p,\Phi_{p}(\eta
^N_{p-1})
)(f_n)|^m |
\Fa^{(N)}_{p-1}
\bigr)^{{1/m}}\leq b(2m)^{2} \delta(R^{\Phi
_{p,n}}) ,
\]
from which we readily conclude that
\begin{eqnarray*}
\EE\bigl(|\sqrt{N}J^N_n(f_n)|^m\bigr)^{{1/m}}&=&
N \EE\Biggl(\Biggl|\sum_{p=0}^n\Ra_{p,n}(\eta^N_p,\Phi
_{p}(\eta
^N_{p-1}))(f_n)
\Biggr|^m \Biggr)^{{1/m}}\\
&\leq& b(2m)^{2} \sum_{p=0}^n \delta
(R^{\Phi_{p,n}}) .
\end{eqnarray*}
Notice that
\[
\sqrt{N} I^N_n=\sum_{p=0}^n\sum_{i=1}^N \Xa^{(N,i)}_{p,n}(f_n)
\qquad\mbox{where }
\Xa^{(N,i)}_{p,n}(f_n)=U^{(N,i)}_p\bigl(\Da_{p,n}^{(N)}(f_n)\bigr),
\]
and the random measures $U^{(N,i)}_p$ are given, for any $g_p\in\mbox
{Osc}_{1}(E_p)$, by
\[
U^{(N,i)}_p(g_p)
:=g_p\bigl(\xi^{(N,i)}_p\bigr)-K_{p,\eta_{p-1}^N}(g_p)\bigl(\xi
^{(N,i)}_{p-1}\bigr) .
\]
In the further development of this section, we fix the final time
horizon $n$ and the
the function $f_n\in\operatorname{Osc}_{1}(E_n)$. To clarify the presentation,
we omit the
final time index and the test function $f_n$, and we set, for any $p$
in $[0,n]$,
\[
X^{(N,i)}_p=\Xa^{(N,i)}_{p,n}(f_n),\qquad
S^N_p=\sum_{q=0}^p\sum_{i=1}^N X^{(N,i)}_{q}
\]
and
\[
R^N_p:=N \sum_{k=0}^p\Ra_{q,n}(\eta^N_q,\Phi_{q}(\eta
^N_{q-1})) .
\]
At the final time horizon, we have
\[
p=n\quad\Longrightarrow\quad S^N_n=\sqrt{N} I^N_n \quad\mbox{and}\quad
R^N_n=\sqrt{N}J^N_n .
\]

By construction, these variables
form a triangular array of $\Fa^N_{p-1}$-condition\-ally independent
random variables and
\[
\EE\bigl(\bigl(X_{p}^{(N,i)}\bigr)^2| \Fa^{N}_{p-1} \bigr)=0 .
\]
In addition, we readily check the following almost sure estimates:
\[
\bigl|X^{(N,i)}_{p}\bigr|\leq\beta(\Da\Phi_{p,n}
)
\quad\mbox{and}\quad
\EE\bigl(\bigl(X_{p}^{(N,i)}\bigr)^2| \Fa^{N}_{p-1}
\bigr)^{{1/2}}\leq\sigma_p
\beta(\Da\Phi_{p,n})\vadjust{\goodbreak}
\]
for any $ 0\leq p\leq n$. The proof of the theorem is now a direct
consequence of Lem\-ma~\ref{lemm}.

\begin{appendix}\label{app}
\section*{Appendix}
\subsection{A first-order composition lemma}\label{secap1}

\setcounter{theo}{0}
\begin{lem} \label{derivform}
For any pair of mappings
$\Phi_1\in\Upsilon(E_0,E_1)$ and $\Phi_2\in\Upsilon(E_1,E_2)$ the
composition
mapping $(\Phi_2\circ\Phi_1)\in\Upsilon(E_0,E_2)$ and we
have the first-order
derivation-type formula
%
%
\setcounter{equation}{0}
\begin{equation}\label{eqcompo}
\Da_{\eta}(\Phi_2\circ\Phi_1)= \Da_{\eta}\Phi
_1 \Da_{\Phi
_1(\eta)}\Phi_2 .
\end{equation}
\end{lem}
\begin{pf}
To check this property, we first observe that
under this condition, we clearly have the
Lipschitz property,
\renewcommand{\theequation}{$\Phi$}
\begin{equation}
\hypertarget{PhiLink}{}
|[ \Phi(\mu)-\Phi(\eta)
](f)|\leq
\int|(\mu-\eta)(h)| T^{\Phi}_{\eta}(f,dh) ,
\end{equation}
for some collection of integral operators $T^{\Phi}_{\eta}$ from $\Ba
(F)$ into the set $\operatorname{Osc}_{1}(E)$ such that
%
%
\setcounter{equation}{1}
\renewcommand{\theequation}{A.\arabic{equation}}
\begin{equation}\label{lip}
\sup_{\eta\in\Pa(E)} \int\operatorname{osc}(h) T^{\Phi}_{\eta
}(f,dh)\leq
\operatorname{osc}(f) \delta(T^{\Phi})
\end{equation}
for some finite constant $\delta(T^{\Phi})<\infty$. Using
this property, we easily check that~(\ref{eqcompo}) is met with
\[
\beta\bigl( \Da(\Phi_2\circ\Phi_1)\bigr)\leq\beta ( \Da \Phi_2) \beta(
\Da\Phi_1)
\]
and
\[
\delta(R^{\Phi_2\circ\Phi_1})\leq\delta
(T^{\Phi
_1})+\delta(T^{\Phi_1})^2 \delta(R^{\Phi
_2}
) .
\]
This ends the proof of the lemma.
\end{pf}

We also mention that for any pair of mappings $\Phi_1\dvtx\eta\in\Pa
(E_0)\mapsto\Phi_1\in\Pa(E_1)$
and $\Phi_2\dvtx\eta\in\Pa(E_1)\mapsto\Phi_1\in\Pa(E_2)$, the
composition mapping
$\Phi=\Phi_2\circ\Phi_1$ satisfies condition \hyperlink{PhiLink}{$(\Phi)$} as soon as this
condition is met for each mapping. In this case, we also notice that
\[
\delta(T^{\Phi_2\circ\Phi_1})\leq\delta
(T^{\Phi
_2})\times\delta(T^{\Phi_1}) .
\]
Suppose we are given a
mapping $\Phi$ defined in terms of a nonlinear transport formula
\[
\Phi(\eta)=\eta K_{\eta} ,
\]
with a collection of Markov transitions $K_{\eta}$ from a measurable
space $(E,\Ea)$
into another $(F,\Fa)$ satisfying condition \hyperlink{lipsKlink}{$(K)$}. Using the decomposition
\[
\Phi(\mu)-\Phi(\eta)=[\eta-\mu]K_{\eta}+\mu
[K_{\mu
}-K_{\eta}] ,\vadjust{\goodbreak}
\]
we readily check that
\[
\mbox{\hyperlink{lipsKlink}{$(K)$}} \Longrightarrow\mbox{\hyperlink{PhiLink}{$(\Phi)$}} \qquad\mbox{with }
T^{\Phi}_{\eta}(f,dh)=\delta_{K_{\eta}(f)}(dh)+T^K_{\eta
}(f,dh) .
\]

\subsection{\texorpdfstring{Proof of Theorem \protect\ref{theocentralfluc}}{Proof of Theorem 5.2}}
\label{secappfluc}

Let $\Fa^N=\{\Fa^N_n ; n\geq0\}$ be the natural filtration
associated with the
$N$-particle system
$\xi^{(N)}_n$. The first class of
martingales that arises naturally in our context is the
$\RR^d$-valued and $\Fa^N$-martingale
$M^{N}_n(f)$ defined by
%
%
\begin{equation}\label{Lamart}
M^{N}_n(f)=\sum_{p=0}^n [\et^N_p(f_p)-\Phi_p(\et
^N_{p-1})(f_p)] ,
\end{equation}
where $ f_p\dvtx x_p\in E_p\mapsto
f_p(x_p)=(f_p^u(x_p))_{u=1,\ldots,d}\in\RR^d $ is a $d$-dimensional and
bounded measurable function. By direct inspection, we see that the
$v$th component of the martingale
$M^{N}_n(f)=(M^{N}_n(f^u))_{u=1,\ldots,d}$ is the $d$-dimensional and
$F^N$-martingale defined for any $u=1,\ldots,d$ by the formula
\begin{eqnarray*}
M^{N}_n(f^u) &=& \sum_{p=0}^n [\et^N_p(f_p^u)-\Phi_p(\et ^N_{p-1})(f_p^u)] \\
&=& \sum_{p=0}^n [\et^N_p(f_p^u)-\et^N_{p-1}K_{p,\eta_{p-1}^N}(f_p^u) ],
\end{eqnarray*}
with the usual\vspace*{-2pt} convention
$K_{0,\et^N_{-1}}=\eta_0=\Phi_0(\et^N_{-1})$ for $p=0$. The idea of the
proof consists of using the CLT for triangular arrays of $\RR^d$-valued
random variables (\cite{JaShli}, Theorem 3.33, page 437). We first
rewrite the martingale $\sqrt{N} M^{N}_n(f)$ in the following form:
\[
\sqrt{N} M^{N}_n(f)=\sum_{i=1}^N\sum_{p=0}^n \frac{1}{\sqrt {N}} \bigl(
f_p\bigl(\xi^{(N,i)}_p\bigr) -K_{p,\eta^N_{p-1}}(f_p)\bigl(\xi^{(N,i)}_{p-1}\bigr) \bigr) .
\]
This readily yields $ \sqrt{N} M^{N}_n(f)=\sum_{k=1}^{(n+1)N}
U^{N}_k(f) $ where for any $1\leq k\leq(n+1)N$ with $k=pN+i$ for some
$i=1,\ldots,N$ and $p=0,\ldots,n$
\[
U^{N}_k(f)=\frac{1}{\sqrt{N}}\bigl( f_p\bigl(\xi^{(N,i)}_p\bigr) -
K_{p,\eta^N_{p-1}}(f_p)\bigl(\xi^{(N,i)}_{p-1}\bigr) \bigr) .
\]
We further denote by $\Ga_k^N$ the $\sigma$-algebra generated by the
random variables $\xi^{j}_p$ for any pair index $(j,p)$ such that
$pN+j\leq k$. It can be checked that, for any $1\leq u<v\leq d$ and for
any $1\leq k\leq(n+1)N$ with $k=pN+i$ for some $i=1,\ldots,N$ and
$p=0,\ldots,n$, we have $ \EE(U^{N}_k(f^u) | \Ga^N_{k-1})=0 $ and
\begin{eqnarray*}
&&\EE(U^{N}_k(f^u)U^{N}_k(f^v) | \Ga^N_{k-1})\\
&&\qquad=\frac{1}{N} K_{p,\eta^N_{p-1}}[(f_p^u -
K_{p,\eta^N_{p-1}}f_p^u) (f_p^v -
K_{p,\eta^N_{p-1}}f_p^v)] \bigl(X^{(N,i)}_{p-1}\bigr) .
\end{eqnarray*}
This also yields that
\begin{eqnarray*}
&&\sum_{k=pN+1}^{pN+N}
\EE(U^{N}_k(f^u)U^{N}_k(f^v) | \Fa^N_{k-1})
\\
&&\qquad=\eta^N_{p-1}\bigl[ K_{p,\eta^N_{p-1}}[(f_p^u -
K_{p,\eta^N_{p-1}}f_p^u) (f_p^v -
K_{p,\eta^N_{p-1}}f_p^v) ]\bigr] .
\end{eqnarray*}
Our aim is now to describe the limiting behavior of the martingale
$\sqrt{N} M^{N}_n(f)$ in terms of the process
$X^{N}_t(f)\stackrel{\mathrm{def}.}{=}\sum_{k=1}^{[Nt]+N} U^{N}_k(f)$.
By the definition of the particle model associated with a given mapping
$\Phi_n$,
and using the fact that
$[\frac{[Nt]}{N}]=[t]$, one gets that for any $1\leq
u,v\leq d$
\begin{eqnarray*}
&&\sum_{k=1}^{[Nt]+N} E(U^{N}_k(f^u) U^{N}_k(f^v)
|\Fa_{k-1}^N )
\\
&&\qquad=C^{N}_{[t]}(f^u,f^v)+\frac{[Nt]-N[t]}{N}
\bigl(C^{N}_{[t]+1}(f^u,f^v) -C^{N}_{[t]}
(f^u,f^v)\bigr),
\end{eqnarray*}
where, for any $n\geq0$ and $1\leq u,v\leq d$,
\[
C^{N}_{n}(f^u,f^v) = \sum_{p=0}^n \eta^N_{p-1}\bigl[
K_{p,\eta^N_{p-1}} \bigl( ( f^u_p- K_{p,\eta^N_{p-1}}f^u_p
)
( f^v_p- K_{p,\eta^N_{p-1}}f^v_p ) \bigr)\bigr] .
\]
Under our regularity conditions on the McKean transitions, this implies
that for any $1\leq i,j\leq d$,
\[
\sum_{k=1}^{[Nt]+N} E(U^{N}_k(f^u) U^{N}_k(f^v)
|\Fa_{k-1}^N )
\mathop{\hbox to 1cm{\rightarrowfill}}^P_{N\rightarrow\infty}
C_t(f^u,f^v) ,
\]
with
\[
C_{n}(f^u,f^v)
=\sum_{p=0}^n \eta_{p-1}\bigl[K_{p,\eta_{p-1}}
\bigl( ( f^u_p-K_{p,\eta_{p-1}}f^u_p )
( f^v_p-K_{p,\eta_{p-1}}f^v_p ) \bigr)\bigr]
\]
and, for any $t\in\RR_+$,
\[
C_t(f^u,f^v) = C_{[t]}(f^u,f^v)+\{t\}
\bigl(C_{[t]+1}(f^u,f^v) -C_{[t]}(f^u,f^v)\bigr).
\]
Since $\|U^{N}_k(f)\|\leq\frac{2}{\sqrt{N}}
(\bigvee_{p\leq n}\|f_p\|)$, for any $1\leq k\leq[Nt]+N$,
the conditional Lindeberg condition is
clearly satisfied, and therefore
one concludes that the $\RR^d$-valued martingale $\{X^{N}_t(f) ;
t\in\RR_+\}$ converges in law to a continuous Gaussian martingale
$\{X_t(f) ; t\in\RR_+\}$ such that, for any $1\leq u,v\leq d$
and $t\in\RR_+$, $
\langle X(f^u), X(f^v) \rangle_t=C_t(f^u,f^v)
$.
Recalling that $X^{N}_{[t]}(f)=\sqrt{N} M^{N}_{[t]}(f)$, we conclude
that the
$\RR^d$-valued and $\Fa^N$-martingale
$\sqrt{N} M^{N}_n(f)$ converges in law to an\vadjust{\goodbreak}
$\RR^d$-valued and Gaussian martingale
$M_n(f)=(M_n(f^u))_{u=1,\ldots,d}$ such that for any $n\geq0$ and
$1\leq u,v\leq d$
\[
\langle M(f^u), M(f^v)
\rangle_n =\sum_{p=0}^n \eta_{p-1}\bigl[K_{p,\eta_{p-1}}
\bigl((f^u_p-K_{p,\eta_{p-1}}f^u_p)(f^v
_p-K_{p,\eta_{p-1}}f^v_p)\bigr)\bigr] ,
\]
with the convention $K_{0,\eta_{-1}}=\eta_0$ for $p=0$.

To take the final step, we let $(\varphi_n)_{n\geq0}$ be a sequence of
bounded measurable functions,
respectively, in $\Ba(E_n)^{d_n}$. We associate with $\varphi=(\varphi_n)_n$
the sequence of
functions $f=(f_p)_{0\leq p\leq n}$ defined for any $0\leq p\leq n$ by
the following formula:
\[
f_p=(f^{u}_p)_{u=0,\ldots,n}=(0,\ldots,0,\varphi_p,0,\ldots,0)\in
\Ba
(E_p)^{d_0+\cdots+d_p+\cdots+d_n}.
\]
In the above display, $0$ stands for the null function in $\Ba
(E_p)^{d_q}$ (for $q\not=p$).
By construction, we have, $f^u_u=\varphi_u$ and for any $0\leq u\leq
n$, we have that
\begin{eqnarray*}
f^{u}&=&(f_p^{u})_{0\leq p\leq n}\\
&=&(0,\ldots,0,\varphi_u,0,\ldots,0)\\
&&{}\in
\Ba
(E_0)^{d_0}\times\cdots\times\Ba(E_u)^{d_u}\times\cdots\times\Ba(E_n)^{d_n}
\end{eqnarray*}
so that
\[
\sqrt{N} M^N_n(f^{u})=
\sqrt{N} [\eta_u^N(\varphi_u)-\eta_{u-1}^NK_{u,\eta
^N_{u-1}}(\varphi_u)]=
V^N_u(\varphi_u)
\]
and therefore
\[
\sqrt{N} M^N_n(f):=\bigl(\sqrt{N} M^N_n(f^{u})\bigr)_{0\leq u\leq n}
=(V^N_u(\varphi_u))_{0\leq u\leq n}:=\Va^N_n(\varphi) .
\]
We conclude that $\Va^N_n(\varphi)$
converges in law to an $(n+1)$-dimensional and centered Gaussian
random field
$
\Va_n(\varphi)=(V_u(\varphi_u))_{0\leq u\leq n}
$
with, for any $0\leq u,v\leq n$,
\begin{eqnarray*}
&&\EE(V_u(\varphi_u^1)V_v(\varphi_v^2))
\\
&&\qquad=1_{u}(v)
\eta_{u-1}[K_{u,\eta_{u-1}}
(\varphi_u^1-K_{u,\eta_{u-1}}\varphi_u^1)
K_{u,\eta_{u-1}}
(\varphi_u^2-K_{u,\eta_{u-1}}\varphi_u^2)] .
\end{eqnarray*}
This ends the proof of the theorem.

\subsection{Feynman--Kac semigroups}
\label{FeynKac}

In the context of Feynman--Kac flows (\ref{fkflows}) discussed in the
\hyperref[intro]{Introduction},
the semigroup $\Phi_{p,n}$ is given by the following formula:
\[
\eta_n(f)=\frac{\eta_p(Q_{p,n}(f))}{\eta_p(Q_{p,n}(1))} \qquad
\mbox{with } Q_{p,n}=Q_{p+1},\ldots, Q_{n-1}Q_n .
\]
For $p=n$, we use the convention $Q_{n,n}=\mathrm{Id}$, the identity operator.
Also observe that
\[
[\Phi_{p,n}(\mu)-\Phi_{p,n}(\eta)](f)=\frac{1}{\mu(G_{p,n,\eta
})} (\mu
-\eta)D_{\eta}\Phi_{p,n}(f) ,
\]
with the first-order operator
\[
D_{\eta}\Phi_{p,n}(f):=G_{p,n,\eta} P_{p,n}\bigl(f-\Phi_{p,n}(\eta
)(f)\bigr) .\vadjust{\goodbreak}
\]
In the above display $G_{p,n,\eta}$ and $P_{p,n}$ stand for the
potential function
and the Markov operator given by
\[
G_{p,n,\eta}:={Q_{p,n}(1)}/{\eta(Q_{p,n}(1))}
\quad\mbox{and}\quad
P_{p,n}(f)={Q_{p,n}(f)}/{Q_{p,n}(1)} .
\]
It is now easy to check that
\[
\Ra^{\Phi_{p,n}}(\mu,\eta)(f):=-\frac{1}{\mu(G_{p,n,\eta})} [\mu
-\eta
]^{\otimes2}\bigl(G_{p,n,\eta}\otimes D_{p,n,\eta}(f)\bigr) .
\]
Using the fact that
\[
D_{\eta}\Phi_{p,n}(f)(x)=G_{p,n,\eta}(x) \int[
P_{p,n}(f)(x)-P_{p,n}(f)(y)] G_{p,n,\eta}(y) \eta(dy) ,
\]
we find that
\[
\forall f\in\operatorname{Osc}_1(E_n)\qquad \Vert D_{\eta}\Phi
_{p,n}(f)\Vert\leq q_{p,n} \beta(P_{p,n})
\]
with
\[
q_{p,n}=\sup_{x,y}\frac{Q_{p,n}(1)(x)}{Q_{p,n}(1)(y)} .
\]
This implies that
\[
\beta(D\Phi_{p,n})\leq2 q_{p,n} \beta(P_{p,n}) .
\]
Finally, we observe that
\[
| \Ra^{\Phi_{p,n}}(\mu,\eta)(f)|\leq
(2 q_{p,n}^2\beta
(D_{p,n})) \biggl|[\mu-\eta]^{\otimes2}\biggl(\frac
{G_{p,n,\eta
}}{2q_{p,n}}\otimes\frac{D_{p,n,\eta}(f)}{\beta(D_{p,n})}\biggr)
\biggr|
\]
from which one concludes that
\[
\delta(R^{\Phi_{p,n}})\leq2 q_{p,n}^2 \beta(D\Phi_{p,n})\leq
4 q_{p,n}^3 \beta(P_{p,n}) .
\]

We end this section with the analysis of these quantities for
the time homogeneous models discussed in~(\ref{fkhom1}) and
(\ref{fkhom2}). Under the condition \hyperlink{fkhom1link}{$(M)_m$} we have for any $n\geq
m\geq1$ and
$p\geq1$,
%
%
\begin{equation} \label{prodd1}%
q_{p,p+n}\leq\delta_{m}/\varepsilon_m \quad\mbox{and}\quad
\beta(P_{p,p+n})\leq( 1-\varepsilon_{m}^{2}/\delta
_{m-1})^{\lfloor n/m\rfloor} .
\end{equation}
The proof of these estimates relies on semigroup techniques (see \cite
{fk}, Chapter~4, for details). Several contraction inequalities can be
deduced from these results, given below.

For any $k\geq0$ and for $l=1,2$,
%
%
\begin{equation}\label{sumcv}
\sum_{p=0}^{n}q_{p,n}^{k} \beta(P_{p,n})^l\leq\varpi_{k,l}(m):=
\frac{m (\delta_{m}/\varepsilon_m)^k}{1-((
1-\varepsilon
_{m}^{2}/\delta_{m-1})
)^l } .
\end{equation}
Notice that
\[
\varpi_{k,l}(m)\leq m \delta_{m-1} \frac{ \delta_{m}^k/\varepsilon
_m^{k+2}}{(2-(\varepsilon_m^2/\delta_{m-1}))^{l-1}}\leq m \delta _{m-1}
\delta_{m}^k/\varepsilon_m^{k+2} ,\vadjust{\goodbreak}
\]
and that
\[
r_n\leq 4 \varpi_{3,1}(m)\quad \mbox{and}\quad b^{\star}_n\leq2
\delta_{m}/\varepsilon_m
\]
as well as
\[
\overline{\sigma}^2_n\leq 4 \varpi_{2,2}(m) \sigma^2
\quad\mbox{and}\quad \beta_n^2\leq 4 \varpi_{2,2}(m) \qquad\mbox{with }
\sigma^2:=\sup_{n\geq1}\sigma^2_n (\mbox{$\leq$} 1) .
\]

\subsection{McKean mean field model of gases}
\label{McKeangas}

We consider McKean-type models of gases~(\ref{nuvelos}) presented in
Section~\ref{gases}. To simplify the presentation, we consider time
homogeneous models, and we supress the time index. In this notation, we
find that
\[
[K_{\eta}-K_{\mu}](f)(x)=\int\nu(ds) [\eta-\mu
]
(a(s,\point)) M(f)(s,x) .
\]
Observe that
\[
[\eta-\mu](K_{\eta}-K_{\mu})(f)(x)=\int\nu(ds) [\eta-\mu]
(a(s,\point))
[\eta-\mu](M(f)(s,\point)) .
\]
Using the decomposition
%
%
\begin{equation}\label{decompform}\quad
\Phi(\eta)-\Phi(\mu)=(\eta-\mu)K_{\mu}+\mu(K_{\eta}-K_{\mu
})+[\eta-\mu
](K_{\eta}-K_{\mu})
\end{equation}
we readily check that $\Phi\in\Upsilon(E,E)$ with the first-order operator
\begin{eqnarray*}
D_{\mu}\Phi(f)(x)&=&[K_{\mu}(f)(x)-\Phi(\mu)(f)]\\
&&{}+\int \nu (ds) [a(s,x)-\mu(a(s,\point))] \mu( M(f)(s,\point))
\end{eqnarray*}
and the second-order remainder measure
\[
\Ra^{\Phi}(\mu,\eta)(f)=\int[\eta-\mu]^{\otimes2}
(g_s) \nu
(ds) \qquad
\mbox{with }
g_s=a(s,\point)\otimes M(f)(s,\point) .
\]
In this situation, we notice that
\[
\beta(\Da\Phi)\leq\beta(M)\biggl[1+\int\nu(ds)
\operatorname{osc}(a(s,\point
))\biggr]
\]
and
\[
\delta(R^{\Phi})\leq\beta(M)\int\nu(ds) \operatorname{osc}(a(s,\point)) .
\]

\subsection{Gaussian semigroups}\label{aGauss}

To simplify the presentation, we only discuss time homogenous and
one-dimensional models.
We consider the one-dimensional gaussian transitions on $E=\RR$
defined below:
\[
K_{\eta}(x,dy)=\frac{1}{\sqrt{2\pi}} \exp{\biggl\{-\frac
{1}{2}
\bigl(y-d(x,\eta)\bigr)^2\biggr\}} \,dy
\]
with some linear drift function $d_n$ of the form $d(x,\eta)=a(x)+\eta
(b) c(x)$, with some measurable (and nonnecessarily bounded) function
$a$, and some\vadjust{\goodbreak} pair of functions
$b$ and $c\in\Ba(\RR)$.
We use the decomposition
\begin{eqnarray*}
[K_{\eta}-K_{\mu}](f)(x)&=&\int K_{\mu}(x,dy) \Theta
(\Delta_{\mu,\eta}(x,y)) f(y)\\
&&{}+
\int K_{\mu}(x,dy) \Delta_{\mu,\eta}(x,y) f(y)
\end{eqnarray*}
with $\Theta(u)=e^u-1-u$ and the function $\Delta_{\mu,\eta}(x,y)$
defined by
\[
\Delta_{\mu,\eta}(x,y)=\log{\frac{dK_{\eta}(x,\point)}{dK_{\mu
}(x,\point)}(y)}=\Delta^{(1)}_{\mu,\eta}(x,y)+\Delta^{(2)}_{\mu
,\eta}(x,y)
\]
with
\begin{eqnarray*}
\Delta^{(1)}_{\mu,\eta}(x,y)&:=&
[d(x,\eta)-d(x,\mu)]
[y-d(x,\mu)]=c(x) (\eta-\mu)(b) [y-d(x,\mu
)],\\
\Delta^{(2)}_{\mu,\eta}(x,y)&:=&-\tfrac{1}{2}
[d(x,\eta)-d(x,\mu)]^2=-\tfrac{1}{2} c(x)^2 [(\eta-\mu
)(b)]^2 .
\end{eqnarray*}
Under our assumptions on the drift function $d$, we have
\[
\bigl|\Delta^{(1)}_{\mu,\eta}(x,y)\bigr|\leq\|c\| \operatorname{osc}(b) |y-d(x,\mu)|
\quad\mbox{and}\quad
\bigl|\Delta^{(2)}_{\mu,\eta}(x,y)\bigr|\leq\|c\|^2 \operatorname{osc}(b)^2/2 .
\]
Using the fact that $|\Theta(u)|\leq e^{|u|}u^2/2$, after
some elementary manipulations
we prove that
\[
\sup_{x\in\RR}\biggl|[K_{\eta}-K_{\mu}](f)(x)-\int
K_{\mu
}(x,dy) \Delta^{(1)}_{\mu,\eta}(x,y) f(y)\biggr|\leq C [(\eta-\mu
)(b)]^2 \|f\|
\]
with some finite constant $C<\infty$ whose values only depend on $ \|
c\|
$ and $\operatorname{osc}(b)$. On the other hand, we have
\[
\int(\eta-\mu) (dx)\int K_{\mu}(x,dy) \Delta^{(1)}_{\mu,\eta}(x,y) f(y)
=(\eta-\mu)^{\otimes2}\bigl(b\otimes(K'_{\mu}(f))\bigr)
\]
and
\[
\int\mu(dx)\int K_{\mu}(x,dy) \Delta^{(1)}_{\mu,\eta}(x,y) f(y)
= (\eta-\mu)(b) \mu(K'_{\mu}(f))
\]
with the bounded integral operator $K'_{\mu}$ defined by
\[
K'_{\mu}(f)(x)=c(x) \int K_{\mu}(x,dy) [y-d(x,\mu)
] f(y) .
\]
Using the decomposition~(\ref{decompform})
we prove that
\[
\Phi(\eta)\bigl(f-\Phi(\mu)(f)\bigr)
=(\eta-\mu) D_{\mu}\Phi(f)+\Ra^{\Phi}(\eta,\mu)(f) ,
\]
with the first-order operator
\[
D_{\mu}\Phi(f)=K_{\mu}\bigl(f-\Phi(\mu)(f)\bigr)+b \mu
(K'_{\mu
}(f))
\]
and a second-order remainder term such that
\[
|\Ra^{\Phi}(\eta,\mu)(f)|\leq C' \bigl[\bigl|(\eta
-\mu
)^{\otimes2}
\bigl(b\otimes(K'_{\mu}(f))\bigr)\bigr|+ [(\eta-\mu
)(b)]^2 \operatorname{osc}(f)\bigr]\vadjust{\goodbreak}
\]
with some finite constant $C'<\infty$ whose values only depend on $ \|
c\|$ and $\operatorname{osc}(b)$.
Using the fact that
\[
K'_{\mu}(1)=0 \quad\mbox{and}\quad
\|K'_{\mu}(f)\|=\bigl\|K'_{\mu}\bigl(f-\Phi(\mu)(f)\bigr) \bigr\|\leq
\|c\| \operatorname{osc}(f) ,
\]
we conclude that~(\ref{condxi222}) and~(\ref{lipr}) are met with
$ \delta(R^{\Phi})\leq C'\operatorname{osc}(b)(2\|c\|+\operatorname
{osc}(b))$, and condition~(\ref{firsto})
is satisfied with $\beta(\Da\Phi)\leq1+\|c\| \operatorname{osc}(b)$.

\subsection{Legendre transform and convex analysis}
\label{sconvex}
We associate with any increasing and convex function $L\dvtx t\in
\operatorname{Dom}(L)\mapsto L(t)\in\RR_+$ defined in some
domain $\operatorname{Dom}(L)\subset\RR_+$, with $L(0)=0$, the
Legendre--Fenchel transform
$L^{\star}$ defined by the variational formula
\[
\forall\lambda\geq0\qquad
L^{\star}(\lambda):=\sup_{t\in\operatorname{Dom}(L)}{ \bigl(\lambda
t-L(t)\bigr)} .
\]
Note that $L^{\star}$ is a convex increasing function with $L^{\star
}(0)=0$ and its inverse
$(L^{\star})^{-1}$ is a concave increasing function [with
$(L^{\star})^{-1}(0)=0$].

For instance, the Legendre--Fenchel transforms $(\alpha_0^{\star
},\alpha
^{\star}_1)$ of the pair of
convex nonnegative functions $(\alpha_0,\alpha_1)$ given below:
\[
\forall t\in[0,1/2[\qquad \alpha_0(t):=-t-\tfrac{1}{2}\log{(1-2t)}
\]
and
\[
\forall t\geq0 \qquad\alpha_1(t):=e^t-1-t
\]
are simply given by
\[
\alpha_0^{\star}(\lambda)=\tfrac{1}{2}\bigl(\lambda-\log
{(1+\lambda
)}\bigr)\quad
\mbox{and}\quad \alpha_1^{\star}(\lambda)=(1+\lambda
)\log
{(1+\lambda)}-\lambda.
\]

Recall that, for any centered random variable $Y$ with values in
$]{-}\infty, 1]$ such that
$\EE(Y^2) \leq v$, we have
%
%
\begin{equation}\label{ben}
\EE(e^{tY}) \leq\frac{ v e ^t + e^{-vt} }{1+v} \leq
1 + v \alpha_1(t) \leq\exp( v\alpha_1 (t) ) .
\end{equation}
We refer to \cite{ben} for a proof of~(\ref{ben}) and for more
precise results.
For any pair of such functions $(L_1,L_2)$, it is readily checked that
\begin{eqnarray*}
&\displaystyle \forall t\in\operatorname{Dom}(L_2)\qquad L_1(t)\leq L_2(t)
\quad\mbox{and}\quad
\operatorname{Dom}(L_2)\subset\operatorname{Dom}(L_1)&
\\
&\Downarrow&
\\
&\displaystyle L^{\star}_2\leq L^{\star}_1 \quad\mbox{and}\quad
(L^{\star}_1)^{-1}\leq(L^{\star}_2)^{-1} .&
\end{eqnarray*}
For any pair of positive numbers $(u,v)$, We also have that
\begin{eqnarray*}
&\displaystyle \forall t\in v^{-1}\operatorname{Dom}(L_2)\qquad
L_1(t)=u L_2(v t)&
\\
&\Downarrow&
\\
&\displaystyle \forall\lambda\geq0 \qquad L^{\star}_1(\lambda)=u L^{\star}_2
\biggl(\frac
{\lambda}{uv}\biggr)\quad
\mbox{and}\quad \forall x\geq0\qquad
( L^{\star}_1)^{-1} (x)=uv ( L^{\star}_2)^{-1} \biggl(\frac
{x}{u}\biggr)
.&
\end{eqnarray*}
As a simple consequence of the latter results, let us quote the
following property
that will be used later in the further development of Section~\ref{secmvp}:
\[
u\leq\overline{u} \quad\mbox{and}\quad
v\leq\overline{v}\quad\Longrightarrow\quad uv ( L^{\star}_2)^{-1} \biggl(\frac
{x}{u}\biggr)\leq
\overline{u} \,\overline{v} ( L^{\star}_2)^{-1}\biggl(\frac
{x}{\overline
{u}}\biggr) .
\]
Here we want to give upper bounds on the inverse functions of the
Legendre transforms.
Our motivation is due to the following result, which avoids the loss of
a factor $2$
when adding exponential inequalities. Let $A$ and $B$ be centered
random variables
with finite log-Laplace transform, which we denote by $\alpha_A$ and
$\alpha_B$, in a
neighborhood of $0$. Then, denoting by $\alpha_{A+B}$ the log-Laplace
transform of $A+B$,
%
%
\begin{equation}\label{breta}
(\alpha_{A+B}^\star)^{-1} (t) \leq(\alpha_A^\star)^{-1} (t) +
(\alpha
_B^\star)^{-1} (t)
\end{equation}
for any positive $t$ (see \cite{rio}, Lemma 2.1).

In order to obtain analytic approximations of these inverse functions,
one can use
the Newton algorithm: let
\[
F (z) = z + \frac{x - \alpha^* (z)}{(\alpha^*)'(z)},
\]
and define the sequence $(z_n)$ by $z_n = F (z_{n-1})$. From the
properties of
the Legendre--Fenchel transform, we also have that
%
%
\begin{equation}\label{DefF}
F(z) = \biggl( \frac{\alpha((\alpha')^{-1} (z))+x}{(\alpha')^{-1} (z)}
\biggr) .
\end{equation}
Now recall the variational formulation of the inverse of the
Legendre--Fenchel transform,
%
%
\begin{equation}\label{invfench}
(\alpha^{\star})^{-1}(x)=
\inf_{t > 0} t^{-1} \bigl(\alpha(t)+x\bigr) ,
\end{equation}
valid for any $x \geq0$ (see \cite{rio2}, page 159 for a proof of this formula).
From this formula, assuming that $\alpha'' (0) >0$ and setting $z=
\alpha'(t)$,
we get that
%
%
\begin{equation}\label{invfench2}
(\alpha^{\star})^{-1}(x)=
\inf_{z\in\alpha' ( \operatorname{Dom} (\alpha))} F(z) .
\end{equation}
Let then $f(z) = \alpha((\alpha')^{-1} (z))+x$ and $g(z) = (\alpha
')^{-1} (z)$.
From the strict convexity of $\alpha$, the function $t \rightarrow
t^{-1} ((\alpha(t)+x)$
a unique minimum $t_x$ and is decreasing with negative derivative
for $t <t_x$,
increasing with positive derivative for $t > t_x$. It follows that
$f/g$ has a
unique critical point $z(x)$, which is the unique global strict minimum
of $F$
and the unique fixed point of $F$. Furthermore $z(x) = (\alpha^\star
)^{-1} (x)$.

Let $z_0>0$ be in the interior of the image by $\alpha'$ of the domain
of $\alpha$.
If $z_0 > z(x)$, then $(z_n)$ is a decreasing sequence of numbers
bounded from below by $z(x)$.
Hence $(z_n)$ decreases to $z(x)$ as $n$ tends to $\infty$. If $z_0< z(x)$
and $F(z_0)$ belongs to the interior of $\alpha' (\operatorname{Dom}
(\alpha) )$,
then $z_1 > z(x)$ and $(z_n)_{n>0}$ is decreasing to $z(x)$.\vadjust{\goodbreak}

We now recall the convergence properties of the Newton algorithm.
Assume that $z_0 > z(x)$ and let $A$ be a positive real such that
$F'' (z) \leq2A$ for any $z$ in $[z(x) , z_0]$. Then, by the Taylor
formula at order $2$,
%
%
\begin{equation}
\label{superG}
0 \leq z_n-z(x) \leq A^{2^n -1} \bigl(z_0 -z(x)\bigr)^{(2^n)} ,
\end{equation}
which provides a supergeometric rate of convergence if $A ( z_0 -z(x) )
< 1$.

Since $F$ depends on $x$, $A$ is a function of $x$. In order to get estimates
of the rate of convergence of $z_n$ to $z(x)$ for small values of $x$,
we now assume that
$\alpha'$ is convex. We will prove that
%
%
\begin{equation} \label{Fseconde}
A := \frac{1}{2} \sup_{z\geq z(x) } F'' (z) \leq\frac{(\alpha
^\star
)^{-1} (x)}{2x\alpha''(0)} .
\end{equation}
To prove~(\ref{Fseconde}),
we start by computing $F'' = (f/g)''$. Since $f'= zg'$,
\[
(f/g)' = g' (zg-f) g^{-2} .
\]
Now $(zg-f)' = g + (zg' -f') = g$. It follows that
\[
(f/g)'' = g' g^{-1} + (zg-f) ( g''g^{-2} - 2 g'^2 g^{-3} ) .
\]
Next, for $z \geq z(x)$, $zg(z) - f(z) \geq0$, so that
\[
(f/g)'' (z) \leq g' g^{-1} + (zg-f) g''g^{-2} .
\]
Under the additional assumption that $\alpha'$ is convex, the inverse function
$(\alpha')^{-1} =g$ is concave, so that $g'' \leq0$. In that case, for
$z\geq z(x)$,
\[
(f/g)'' (z) \leq g' (z) / g(z) = (\log g)' (z) .
\]
Now $\log g$ is the inverse function of $\psi(t) = \alpha' (e^t)$.
From the
properties of $\alpha'$, the function $\psi$ is convex, so that $\log
g$ is
concave. Hence $(\log g)'$ is nonincreasing, which implies that
\[
F'' (z) \leq g' (z(x)) / g(z(x)) = z(x) g' (z(x) )/ f(z(x)) \qquad\mbox{for any } z \geq z(x) .
\]
Since $f(z) \geq x$ and $g' ( z(x) ) \leq g'(0) = 1/\alpha'' (0)$, we get
(\ref{Fseconde}), noticing that $z(x) = F(z(x)) = (\alpha^\star
)^{-1} (x)$.

We now apply these results to the functions $\alpha_0$ and $\alpha_1$.
Using the fact~that
\[
\frac{t^2}{2} \leq\alpha_1(t):=e^t-1-t \leq\overline{\alpha
}_1(t):=\frac{t^2}{2(1-t/3)}
\]
for every $t\in[0,3[$, and applying (B.5), page 153 in \cite{rio2}, we
get that
\[
\sqrt{2x} \leq(\alpha^{\star}_1)^{-1}(x)\leq
(\overline
{\alpha}^{\star}_1)^{-1}(x)
= \sqrt{2x} + (x/3) .
\]
Also, by the second part of Theorem B.2 in \cite{rio2}, the function
$\overline{\alpha}^{\star}_1$, which is the
inverse function of the above function, satisfies
%
%
\begin{equation}\label{compBern}
\overline{\alpha}^{\star}_1 (t) \geq\frac{t^2}{2 (1 + (t/3) ) } ,\vadjust{\goodbreak}
\end{equation}
which is the usual bound in the Bernstein inequality.
Now $z = e^t - 1$, and consequently $t = \log(1+z)$ and
\[
F(z) = \frac{x + z - \log(1+z)}{\log(1+z)} .
\]
Set $z_0 = \sqrt{2x} + (x/3)$. Then $z_0 > z(x)$. Hence $z(x) < z_1 < z_0$
(here $z_1 = F(z_0)$). So
%
%
\begin{eqnarray}
\label{majstar1}
(\alpha^{\star}_1)^{-1}(x) &\leq& z_1 :=
\frac{ \sqrt{2x} + (4x/3) - \log(1 +(x/3) + \sqrt{2x})}{\log(1 +(x/3)
+ \sqrt{2x})}\nonumber\\[-8pt]\\[-8pt]
&\leq&(x/3) + \sqrt{2x} .\nonumber
\end{eqnarray}
Furthermore, from~(\ref{superG}) and~(\ref{Fseconde}) and the fact that
$z_0 - z(x) \leq x/3$,
\[
0 \leq z_1 - (\alpha^{\star}_1)^{-1}(x) \leq\frac{x}{18}
(\alpha^{\star}_1)^{-1}(x) ,
\]
which ensures that
\[
18 z_1 /(18 + x) \leq(\alpha^{\star}_1)^{-1}(x) \leq z_1.
\]
In the same way, noticing that
\[
t^2/(1-4t/3) \leq\alpha_0 (t) \leq t^2/(1-2t) \qquad\mbox{for any } t
\in
[0,1/2 [ ,
\]
we get
\[
2\sqrt{x} + (4x/3) \leq(\alpha_0^\star)^{-1} (x) \leq2\sqrt{x} + 2x
:= z_0.
\]
By definition of $\alpha_0$, we have $\alpha'_0 (t) = 2t/(1-2t)$.
Let $z = 2t/(1-2t)$. Then $t = z/(2+2z)$, so that
\[
F(z) = \frac{x + \alpha_0 ((\alpha'_0)^{-1} (t) )}{(\alpha'_0)^{-1} (t)}
= 2x + \log( 1 + z) + \frac{2x + \log(1+z) - z}{z} .
\]
Computing $z_1 = F(z_0)$, we get
%
%
\begin{eqnarray}
\label{majstar0}
(\alpha_0^\star)^{-1} (x) &\leq& z_1:= 2x + \log\bigl( 1 + 2x+ 2\sqrt{x} \bigr)\nonumber\\
&&\hspace*{23.5pt}{} +
\frac{\log( 1 + 2x+ 2\sqrt{x} ) - 2\sqrt{x} } { 2x + 2\sqrt{x}}\\
&\leq&
2x + 2\sqrt{x} ,\nonumber
\end{eqnarray}
which improves on the previous upper bound. Furthermore,
from~(\ref{superG}) and~(\ref{Fseconde})
\[
0 \leq z_1 - (\alpha^{\star}_0)^{-1}(x) \leq\frac
{x}{9}
(\alpha^{\star}_0)^{-1}(x) ,
\]
which ensures that
\[
9 z_1 /(9 + x) \leq(\alpha^{\star}_0)^{-1}(x) \leq
z_1 .
\]

\end{appendix}

%

%
\printaddresses


\begin{thebibliography}{32}

\bibitem{Bolley}
%
\begin{barticle}[vtex]
\bauthor{\bsnm{Bolley},~\bfnm{F.}\binits{F.}},
\bauthor{\bsnm{Guillin},~\bfnm{A.}\binits{A.}} \AND
\bauthor{\bsnm{Malrieu},~\bfnm{F.}\binits{F.}}
(\byear{2010}).
\btitle{Trend to equilibrium and particle approximation for
a weakly selfconsistent Vlasov--Fokker--Planck equation}.
\bjournal{M2AN}
\bvolume{44}
\bpages{867--884}.
\end{barticle}%
%
\endbibitem%

\bibitem{Villani}
%
\begin{barticle}[mr]
\bauthor{\bsnm{Bolley},~\bfnm{Fran{\c{c}}ois}\binits{F.}},
\bauthor{\bsnm{Guillin},~\bfnm{Arnaud}\binits{A.}} \AND
\bauthor{\bsnm{Villani},~\bfnm{C{\'e}dric}\binits{C.}}
(\byear{2007}).
\btitle{Quantitative concentration inequalities for empirical measures on
non-compact spaces}.
\bjournal{Probab. Theory Related Fields}
\bvolume{137}
\bpages{541--593}.
\bid{doi={10.1007/s00440-006-0004-7}, mr={2280433}}
\end{barticle}
%
\endbibitem

\bibitem{bossy1}
%
\begin{barticle}[vtex]
\bauthor{\bsnm{Bossy},~\bfnm{Mireille}\binits{M.}} \AND
\bauthor{\bsnm{Talay},~\bfnm{Denis}\binits{D.}}
(\byear{1997}).
\btitle{A stochastic particle method for the {M}c{K}ean--{V}lasov and the
{B}urgers equation}.
\bjournal{Math. Comp.}
\bvolume{66}
\bpages{157--192}.
\bid{doi={10.1090/S0025-5718-97-00776-X}, mr={1370849}}
\end{barticle}
%
\endbibitem

\bibitem{bossy2}
%
\begin{barticle}[mr]
\bauthor{\bsnm{Bossy},~\bfnm{Mireille}\binits{M.}} \AND
\bauthor{\bsnm{Talay},~\bfnm{Denis}\binits{D.}}
(\byear{1996}).
\btitle{Convergence rate for the approximation of the limit law of weakly
interacting particles: Application to the {B}urgers equation}.
\bjournal{Ann. Appl. Probab.}
\bvolume{6}
\bpages{818--861}.
\bid{doi={10.1214/aoap/1034968229}, mr={1410117}}
\end{barticle}
%
\endbibitem

\bibitem{cerou}
%
\begin{bmisc}[vtex]
\bauthor{\bsnm{Cerou},~\bfnm{Fr.}\binits{F.}},
\bauthor{\bsnm{Del~Moral},~\bfnm{P.}\binits{P.}} \AND
\bauthor{\bsnm{Guyader},~\bfnm{A.}\binits{A.}}
(\byear{2010}).
\bhowpublished{A non asymptotic variance theorem for unnormalized
Feynman-Kac particle models. 
\textit{Ann. Inst. H. Poincar\'e}. To appear}.
\end{bmisc}
%
\endbibitem

\bibitem{ben}
%
\begin{barticle}[mr]
\bauthor{\bsnm{Bentkus},~\bfnm{Vidmantas}\binits{V.}}
(\byear{2004}).
\btitle{On {H}oeffding's inequalities}.
\bjournal{Ann. Probab.}
\bvolume{32}
\bpages{1650--1673}.
\bid{doi={10.1214/009117904000000360}, mr={2060313}}
\end{barticle}
%
\endbibitem

\bibitem{chopine}
%
\begin{barticle}[mr]
\bauthor{\bsnm{Chopin},~\bfnm{Nicolas}\binits{N.}}
(\byear{2004}).
\btitle{Central limit theorem for sequential {M}onte {C}arlo methods
and its
application to {B}ayesian inference}.
\bjournal{Ann. Statist.}
\bvolume{32}
\bpages{2385--2411}.
\bid{doi={10.1214/009053604000000698}, mr={2153989}}
\end{barticle}
%
\endbibitem

\bibitem{Cro}
%
\begin{barticle}[mr]
\bauthor{\bsnm{Crooks},~\bfnm{Gavin~E.}\binits{G.~E.}}
(\byear{1998}).
\btitle{Nonequilibrium measurements of free energy differences for
microscopically reversible {M}arkovian systems}.
\bjournal{J. Stat. Phys.}
\bvolume{90}
\bpages{1481--1487}.
\bid{doi={10.1023/A:1023208217925}, mr={1628273}}
\end{barticle}
%
\endbibitem

\bibitem{fk}
%
\begin{bbook}[vtex]
\bauthor{\bsnm{Del~Moral},~\bfnm{Pierre}\binits{P.}}
(\byear{2004}).
\btitle{Feynman--{K}ac Formulae:
Genealogical and Interacting Particle Systems with Applications}.
\bpublisher{Springer}, \baddress{New York}.
\bid{mr={2044973}}
\end{bbook}
%
\endbibitem

\bibitem{ddj}
%
\begin{barticle}[mr]
\bauthor{\bsnm{Del~Moral},~\bfnm{Pierre}\binits{P.}},
\bauthor{\bsnm{Doucet},~\bfnm{Arnaud}\binits{A.}} \AND
\bauthor{\bsnm{Jasra},~\bfnm{Ajay}\binits{A.}}
(\byear{2006}).
\btitle{Sequential {M}onte {C}arlo samplers}.
\bjournal{J. R. Stat. Soc. Ser. B Stat. Methodol.}
\bvolume{68}
\bpages{411--436}.
\bid{doi={10.1111/j.1467-9868.2006.00553.x}, mr={2278333}}
\end{barticle}
%
\endbibitem

\bibitem{Miclo3}
%
\begin{barticle}[vtex]
\bauthor{\bsnm{Del~Moral},~\bfnm{Pierre}\binits{P.}} \AND
\bauthor{\bsnm{Miclo},~\bfnm{L.}\binits{L.}}
(\byear{2003}).
\btitle{Particle approximations of {L}yapunov exponents connected to
{S}chr\"odinger operators and {F}eynman--{K}ac semigroups}.
\bjournal{ESAIM Probab. Stat.}
\bvolume{7}
\bpages{171--208}.
\bid{doi={10.1051/ps:2003001}, mr={1956078}}
\end{barticle}
%
\endbibitem

\bibitem{dmiclo}
%
\begin{bincollection}[vtex]
\bauthor{\bsnm{Del~Moral},~\bfnm{P.}\binits{P.}} \AND
\bauthor{\bsnm{Miclo},~\bfnm{L.}\binits{L.}}
(\byear{2000}).
\btitle{Branching and interacting particle systems approximations of
{F}eynman--{K}ac formulae with applications to non-linear filtering}.
In \bbooktitle{S\'eminaire de {P}robabilit\'es, {XXXIV}}.
\bseries{Lecture Notes in Math.}
\bvolume{1729}
\bpages{1--145}.
\bpublisher{Springer}, \baddress{Berlin}.
\bid{doi={10.1007/BFb0103798}, mr={1768060}}
\end{bincollection}
%
\endbibitem

\bibitem{DDS}
%
\begin{barticle}[vtex]
\bauthor{\bsnm{Del~Moral},~\bfnm{Pierre}\binits{P.}} \AND
\bauthor{\bsnm{Doucet},~\bfnm{Arnaud}\binits{A.}}
(\byear{2004}).
\btitle{Particle motions in absorbing medium with hard and soft obstacles}.
\bjournal{Stoch. Anal. Appl.}
\bvolume{22}
\bpages{1175--1207}.
\bid{doi={10.1081/SAP-200026444}, mr={2089064}}
\end{barticle}
%
\endbibitem

\bibitem{DG}
%
\begin{barticle}[mr]
\bauthor{\bsnm{Del~Moral},~\bfnm{P.}\binits{P.}} \AND
\bauthor{\bsnm{Guionnet},~\bfnm{A.}\binits{A.}}
(\byear{1999}).
\btitle{Central limit theorem for nonlinear filtering and interacting particle
systems}.
\bjournal{Ann. Appl. Probab.}
\bvolume{9}
\bpages{275--297}.
\bid{doi={10.1214/aoap/1029962742}, mr={1687359}}
\end{barticle}
%
\endbibitem

\bibitem{arnaud}
%
\begin{bbook}[vtex]
\bauthor{\bsnm{Doucet},~\bfnm{Arnaud}\binits{A.}},
\bauthor{\bsnm{de Freitas},~\bfnm{Nando}\binits{N.}} \AND
\bauthor{\bsnm{Gordon},~\bfnm{Neil}\binits{N.}}
(\byear{2001}).
\btitle{Sequential {M}onte {C}arlo Methods in Practice}.
\bpublisher{Springer}, \baddress{New York}.
\bid{mr={1847783}}
\end{bbook}
%
\endbibitem

\bibitem{steins}
%
\begin{barticle}[mr]
\bauthor{\bsnm{Steinsaltz},~\bfnm{David}\binits{D.}},
\bauthor{\bsnm{Evans},~\bfnm{Steven~N.}\binits{S.~N.}} \AND
\bauthor{\bsnm{Wachter},~\bfnm{Kenneth~W.}\binits{K.~W.}}
(\byear{2005}).
\btitle{A generalized model of mutation-selection balance with
applications to
aging}.
\bjournal{Adv. in Appl. Math.}
\bvolume{35}
\bpages{16--33}.
\bid{doi={10.1016/j.aam.2004.09.003}, mr={2141503}}
\end{barticle}
%
\endbibitem

\bibitem{harris}
%
\begin{barticle}[vtex]
\bauthor{\bsnm{Harris},~\bfnm{T.~E.}\binits{T.~E.}} \AND
\bauthor{\bsnm{Kahn},~\bfnm{H.}\binits{H.}}
(\byear{1951}).
\btitle{Estimation of particle transmission by random
sampling}.
\bjournal{Natl. Bur. Stand. Appl. Math. Ser.}
\bvolume{12}
\bpages{27--30}.
\end{barticle}
%
\endbibitem

\bibitem{JaShli}
%
\begin{bbook}[mr]
\bauthor{\bsnm{Jacod},~\bfnm{Jean}\binits{J.}} \AND
\bauthor{\bsnm{Shiryaev},~\bfnm{Albert~N.}\binits{A.~N.}}
(\byear{1987}).
\btitle{Limit Theorems for Stochastic Processes}.
\bseries{Grundlehren der Mathematischen Wissenschaften [Fundamental Principles
of Mathematical Sciences]}
\bvolume{288}.
\bpublisher{Springer}, \baddress{Berlin}.
\bid{mr={0959133}}
\end{bbook}
%
\endbibitem

\bibitem{Jar}
%
\begin{barticle}[vtex]
\bauthor{\bsnm{Jarzynski},~\bfnm{C.}\binits{C.}}
(\byear{1997}).
\btitle{Nonequilibrium equality for free energy differences}.
\bjournal{Phys. Rev. Lett.}
\bvolume{78}
\bpages{2690--2693}.
\end{barticle}
%
\endbibitem

\bibitem{Jar2}
%
\begin{barticle}[vtex]
\bauthor{\bsnm{Jarzynski},~\bfnm{C.}\binits{C.}}
(\byear{1997}).
\btitle{Equilibrium free-energy differences from nonequilibrium
measurements: A master-equation approach}.
\bjournal{Phys. Rev. E}
\bvolume{56}
\bpages{5018}.
\end{barticle}
%
\endbibitem

\bibitem{Malrieu2}
%
\begin{barticle}[mr]
\bauthor{\bsnm{Malrieu},~\bfnm{Florent}\binits{F.}}
(\byear{2003}).
\btitle{Convergence to equilibrium for granular media equations and their
{E}uler schemes}.
\bjournal{Ann. Appl. Probab.}
\bvolume{13}
\bpages{540--560}.
\bid{doi={10.1214/aoap/1050689593}, mr={1970276}}
\end{barticle}
%
\endbibitem

\bibitem{Talay}
%
\begin{bincollection}[mr]
\bauthor{\bsnm{Malrieu},~\bfnm{Florent}\binits{F.}} \AND
\bauthor{\bsnm{Talay},~\bfnm{Denis}\binits{D.}}
(\byear{2005}).
\btitle{Concentration inequalities for {E}uler schemes}.
In \bbooktitle{Monte {C}arlo and Quasi-{M}onte {C}arlo Methods 2004}
\bpages{355--371}.
\bpublisher{Springer}, \baddress{Berlin}.
\bid{doi={10.1007/3-540-31186-6_21}, mr={2208718}}
\bptnote{check year}
\end{bincollection}
%
\endbibitem

\bibitem{mckean}
%
\begin{bincollection}[mr]
\bauthor{\bsnm{McKean},~\bfnm{H.~P.}\binits{H.~P.}\bsuffix{~Jr.}}
(\byear{1967}).
\btitle{Propagation of chaos for a class of non-linear parabolic equations}.
In \bbooktitle{Stochastic {D}ifferential {E}quations ({L}ecture
{S}eries in
{D}ifferential {E}quations, {S}ession 7, {C}atholic {U}niv., 1967)}
\bpages{41--57}.
\bpublisher{Air Force Office Sci. Res.}, \baddress{Arlington, VA}.
\bid{mr={0233437}}
\end{bincollection}
%
\endbibitem

\bibitem{sylvie}
%
\begin{bincollection}[vtex]
\bauthor{\bsnm{M{\'e}l{\'e}ard},~\bfnm{Sylvie}\binits{S.}}
(\byear{1996}).
\btitle{Asymptotic behaviour of some interacting particle systems;
{M}c{K}ean--{V}lasov and {B}oltzmann models}.
In \bbooktitle{Probabilistic Models for Nonlinear Partial Differential
Equations ({M}ontecatini {T}erme, 1995)}.
\bseries{Lecture Notes in Math.}
\bvolume{1627}
\bpages{42--95}.
\bpublisher{Springer}, \baddress{Berlin}.
\bid{doi={10.1007/BFb0093177}, mr={1431299}}
\end{bincollection}
%
\endbibitem

\bibitem{rio}
%
\begin{barticle}[mr]
\bauthor{\bsnm{Rio},~\bfnm{Emmanuel}\binits{E.}}
(\byear{1994}).
\btitle{Local invariance principles and their application to density
estimation}.
\bjournal{Probab. Theory Related Fields}
\bvolume{98}
\bpages{21--45}.
\bid{doi={10.1007/BF01311347}, mr={1254823}}
\end{barticle}
%
\endbibitem

\bibitem{rio2}
%
\begin{bbook}[vtex]
\bauthor{\bsnm{Rio},~\bfnm{Emmanuel}\binits{E.}}
(\byear{2000}).
\btitle{Th\'eorie Asymptotique des Processus Al\'eatoires Faiblement
D\'ependants}.
\bseries{Math\'ematiques and Applications (Berlin) [Mathematics and
Applications]}
\bvolume{31}.
\bpublisher{Springer}, \baddress{Berlin}.
\bid{mr={2117923}}
\end{bbook}
%
\endbibitem

\bibitem{rosen}
%
\begin{barticle}[vtex]
\bauthor{\bsnm{Rosenbluth},~\bfnm{M.~N.}\binits{M.~N.}} \AND
\bauthor{\bsnm{Rosenbluth},~\bfnm{A.~W.}\binits{A.~W.}}
(\byear{1955}).
\btitle{Monte-Carlo calculations of the average extension of
macromolecular chains}.
\bjournal{J. Chem. Phys.}
\bvolume{23}
\bpages{356--359}.
\end{barticle}
%
\endbibitem

\bibitem{Rousset2}
%
\begin{barticle}[mr]
\bauthor{\bsnm{Rousset},~\bfnm{Mathias}\binits{M.}}
(\byear{2006}).
\btitle{On the control of an interacting particle estimation of
{S}chr\"odinger
ground states}.
\bjournal{SIAM J. Math. Anal.}
\bvolume{38}
\bpages{824--844 (electronic)}.
\bid{doi={10.1137/050640667}, mr={2262944}}
\end{barticle}
%
\endbibitem

\bibitem{Rubin1}
%
\begin{barticle}[mr]
\bauthor{\bsnm{Rubinstein},~\bfnm{Reuven}\binits{R.}}
(\byear{2009}).
\btitle{The {G}ibbs cloner for combinatorial optimization, counting and
sampling}.
\bjournal{Methodol. Comput. Appl. Probab.}
\bvolume{11}
\bpages{491--549}.
\bid{doi={10.1007/s11009-008-9101-7}, mr={2551567}}
\end{barticle}
%
\endbibitem

\bibitem{Rubin2}
%
\begin{barticle}[vtex]
\bauthor{\bsnm{Rubinstein},~\bfnm{R.~Y.}\binits{R.~Y.}}
(\byear{2010}).
\btitle{Randomized algorithms with splitting: Why the classic
randomized algorithms do not work and how to make them work}.
\bjournal{Methodol. Comput. Appl. Probab.}
\bvolume{12}
\bpages{1--50}.
\end{barticle}
%
\endbibitem

\bibitem{Shiga}
%
\begin{barticle}[mr]
\bauthor{\bsnm{Shiga},~\bfnm{Tokuzo}\binits{T.}} \AND
\bauthor{\bsnm{Tanaka},~\bfnm{Hiroshi}\binits{H.}}
(\byear{1985}).
\btitle{Central limit theorem for a system of {M}arkovian particles
with mean
field interactions}.
\bjournal{Z. Wahrsch. Verw. Gebiete}
\bvolume{69}
\bpages{439--459}.
\bid{doi={10.1007/BF00532743}, mr={0787607}}
\end{barticle}
%
\endbibitem

\bibitem{Sznitman}
%
\begin{bincollection}[vtex]
\bauthor{\bsnm{Sznitman},~\bfnm{Alain-Sol}\binits{A.-S.}}
(\byear{1991}).
\btitle{Topics in propagation of chaos}.
In \bbooktitle{\'{E}cole D'\'{E}t\'e de {P}robabilit\'es de {S}aint--{F}lour
{XIX}---1989}.
\bseries{Lecture Notes in Math.}
\bvolume{1464}
\bpages{165--251}.
\bpublisher{Springer}, \baddress{Berlin}.
\bid{doi={10.1007/BFb0085169}, mr={1108185}}
\end{bincollection}
%
\endbibitem

\end{thebibliography}
\end{document}